\input amstex
\documentstyle{amsppt}
\magnification=1200
\nologo
\NoBlackBoxes
\input xy
\xyoption{all}
\topmatter
\title De Rham and infinitesimal cohomology in Kapranov's model for noncommutative 
algebraic geometry\endtitle

\author Guillermo Corti\~nas$^*$\endauthor

\affil Departamento de Matem\'atica\\
       Facultad de Ciencias Exactas y Naturales\\
       Universidad de Buenos Aires\endaffil 

\address Departamento de Matem\'atica, Ciudad Universitaria Pabell\'on 1,
(1428) Buenos Aires, Argentina\endaddress
\email gcorti\@dm.uba.ar\endemail
\abstract The title refers to the nilcommutative or $NC$-schemes introduced by M. Kapranov
in {\it Noncommutative geometry based on commutator expansions,} J. reine angew. Math {\bf 505} (1998) 
73-118. The latter are noncommutative nilpotent thickenings of commutative schemes. We consider also the 
parallel theory of nil-Poisson or $NP$-schemes, which are nilpotent thickenings of commutative 
schemes in the category of Poisson schemes.
We study several variants of de Rham cohomology for $NC$- and $NP$-schemes. The variants include 
nilcommutative and nil-Poisson versions of the de Rham complex as well as of the 
cohomology of the infinitesimal site introduced by Grothendieck in 
{\it Crystals and the de Rham cohomology of schemes,} Dix expos\'es sur la cohomologie
des sch\'emas, Masson \& Cie, North-Holland (1968) 306-358. It turns out that each of these noncommutative
variants admits a kind of Hodge decomposition which allows one to express the cohomology groups of a 
noncommutative scheme $Y$ as a sum of copies of the usual (de Rham, infinitesimal) cohomology groups of 
the underlying commutative scheme $X$ (Theorems 6.2, 6.5, 6.8).
As a byproduct we obtain new proofs for classical results of Grothendieck (Corollary 6.3) and of 
Feigin-Tsygan (Corollary 6.9) on the relation between de Rham and infinitesimal cohomology and between 
the latter and periodic cyclic homology. 
\endabstract
\thanks (*) Researcher at CONICET. Partially supported by grants 
BID802/OC-AR-PICT 2260 and UBACyT TW79. Part of this research was carried out while
visiting the universities of M\"unster and Bielefeld with a Humboldt fellowship.\endthanks
\rightheadtext{Noncommutative De Rham cohomology}
\leftheadtext{Guillermo Corti\~nas}
\endtopmatter

\define\op{\operatorname}
\define\sg{\operatorname{sg}}
\define\Ass{\op{Ass}}
\define\comm{\op{Comm}}

\define\Hom{\op{Hom}}
\define\cyl{\op{Cyl}}

\define\Poiss{\op{Poiss}}
\define\coli{\operatornamewithlimits{colim}}
\define\holi{\operatornamewithlimits{holim}}

\define\Spec{\op{Spec}}
\define\pro{\op{Pro}}
\define\zar{\op{Zar}}
\define\coker{\op{coker}}

\define\nyl{\overline{\cyl}}
\define\ncl{{{NC}_l}}
\define\nci{{{NC}_\infty}}
\define\ncm{{{NC}_m}}
\define\npl{{{NP}_l}}
\define\npi{{{NP}_\infty}}

\define\calo{\Cal O}
\define\calg{\Cal G}
\define\cals{\Cal S}
\define\calc{\Cal C}
\define\caly{\Cal Y}
\define\calz{\Cal Z}
\define\calw{\Cal W}
\define\calu{\Cal U}
\define\chas{\calc_{\caly}(\cals)}

\define\g{\frak{g}}
\define\N{\frak{N}}
\define\fx{\frak{X}}
\define\fy{\frak{Y}}
\define\fp{\frak{p}}
\define\fq{\frak{q}}
\define\I{\frak{I}}
\define\fib{\twoheadrightarrow}

\define\harrow{\hookrightarrow}
\define\iso{\overset\sim\to\to}
\define\Hy{{\Bbb{H}}}
\document
\head{1. Introduction}\endhead
\bigskip
In this paper we study the de Rham theory of the 
nilcommutative or $NC$-schemes introduced by Kapranov in [Kap]. 
To start let us recall the definitions of $NC$-algebras and schemes
and introduce differential forms for such objects. 
We consider algebras and schemes over a fixed field $k$ of characteristic zero. 
Recall an associative algebra $R$ is {\it nilcommutative of order} $\le l$ or an $\ncl$-algebra if for 
the {\it commutator filtration}
$$
F_0R=R,\qquad F_{n+1}R:=\sum_{p=1}^nF_pRF_{n+1-p}R+\sum_{p=0}^n
<[F_pR,F_{n-p}R]>\tag{1}
$$
we have $F_{l+1}R=0$. 
For every $\ncl$-algebra $R$ there is defined a noncommutative locally ringed space
$\Spec R$. Its underlying topological space is the prime spectrum of the commutative algebra $A=R/F_1R$;
the stalk at $\fp\in \Spec A$ is the \O re localization of $R$ at the inverse image of $\fp$ under the
projection $R\fib A$. Prime spectra of $\ncl$-algebras are called {\it affine $\ncl$-schemes}. 
In general an $\ncl$-scheme is a locally ringed space that can be covered by affine $\ncl$-schemes.
If $Y=(Y,\calo_{Y})$ is an $\ncl$-scheme then 
$Y^{[0]}=(Y,\calo_Y/F_1\calo_Y)$ is a commutative (i.e. usual) scheme. In particular an 
$NC_0$-scheme is just a commutative or $\comm$-scheme. There is a natural notion
of differential forms for $\ncl$-schemes, as follows. 
For $R\in\ncl$ we define its $\ncl$-DGA of differental forms as the quotient
$$
\Omega_{\ncl}R=\frac{\Omega R}{F_{l+1}\Omega R}\tag{2}
$$
of the usual $DGA$ of noncommutative forms ([CQ1]) by the $l+1$-th term of its commutator filtration
(taken in the $DG$ sense). For example if $R$ is commutative and $l=0$, then this is the usual commutative $DGA$
of K\"ahler forms. One checks that $\Omega_{\ncl}$ localizes (in the \O re sense), and thus
defines a sheaf of $\ncl$-DGA's on $\Spec R$ which is quasi-coherent in the sense that each term
$\Omega_{\ncl}^p$ comes from an $R$-bimodule and its (\O re) localizations (see sections 2 and 4 below). 
Thus for every $\ncl$-scheme $Y$ there is defined a quasi-coherent sheaf $\Omega_{\ncl}$ of 
$\ncl$-DGA's. We compute its cohomology in the formally smooth case. Recall
from [Kap] that an $\ncl$-algebra $R$ is formally $\ncl$-smooth if $\hom_{\ncl}(R,\cdot)$ carries surjections with nilpotent kernel into surjections. We call an $\ncl$-scheme formally $\ncl$-smooth if it can
be covered by spectra of formally $\ncl$-smooth algebras. We remark that $Y$ formally $\ncl$-smooth $\Rightarrow$
$Y^{[0]}$ formally $\comm$-smooth. We show (Corollary 6.3) that if $Y$ is formally 
$\ncl$-smooth then for $X=Y^{[0]}$
$$
\Hy^*(Y_{\zar},\Omega_{\ncl}):=\Hy^*(X_{\zar},\Omega_{\ncl}\calo_Y)=
\Hy^*(X_{\zar},\Omega_{\comm})=:H_{dR}^*X\tag{3}
$$ 
is just the usual de Rham cohomology of the underlying commutative scheme. 
In the affine case, because of the quasi-coherence of the Zariski sheaf
$\Omega_{\ncl}$ its hypercohomology is just the cohomology of its global sections (cf. \thetag{27}) 
and we have
$$
H^*(\Omega_{\ncl}R)=\Hy^*(\Spec R,\Omega_{\ncl})=\Hy^*(\Spec A,\Omega_{\comm})=H^*(\Omega_{\comm}A)=
:H_{dR}^*A
$$
This contrasts with the fact that for every $R\in\Ass$ the usual $DGA$ of noncommutative
forms is acyclic, that is
$$
H^n (\Omega R)=\left\{\matrix k & \text{ if }  n=0\\
                              0 &\text { if } n\ne 0\endmatrix\right.
$$
Recall however that if we divide $\Omega R$ by its commutator subspace we get the image of the periodicity
map in cyclic homology ([L, 2.6.7])
$$
H^n(\frac{\Omega R}{[\Omega R,\Omega R]})=s(HC_{n+2}R)\subset HC_nR
$$
which is nontrivial in general. For example if $A$ is smooth commutative then ([L 5.1.12])
$$
H^n(\frac{\Omega A}{[\Omega A,\Omega A]})=\bigoplus_{2m\le n}H_{dR}^{n-2m}A
$$
We show (see 8.3.4 below) that if $R$ is formally $\ncl$-smooth then for $A=R/F_1R$,
$$
H^n(\frac{\Omega_{`\ncl}R}{[\Omega_{\ncl}R,\Omega_{\ncl}R]})=\bigoplus_{m=0}^lH_{dR}^{n+2m}A
$$
Note that $H_{dR}^{n+2m}A$ is the $n$-th cohomology of the complex
$$
\tau_{2m}\Omega_{\comm}A:\frac{\Omega^{2m}_{\comm}A}{d\Omega^{2m-1}_{\comm}A}\to\Omega_{\comm}^{2m+1}A\to
\dots\tag{4}
$$
which has $\Omega^{2m}_{\comm}A/d\Omega^{2m-1}_{\comm}A$ in degree zero. We show (Corollary 6.6) 
that if $Y$ is a formally
$\ncl$-smooth scheme then for $X=Y^{[0]}$
$$
\Hy^n(Y_{\zar},\frac{\Omega_{\ncl}}{[\Omega_{\ncl},\Omega_{\ncl}]})=
\bigoplus_{m=0}^l\Hy^{n+2m}(X_{\zar},\tau_{2m}\Omega_{\comm})
\tag{5}
$$
Here $X=Y^{[0]}$, and $\Omega_{\ncl}/[\Omega_{\ncl},\Omega_{\ncl}]$ and
$\tau_{2m}\Omega_{\comm}$ are the sheafified complexes.
In particular the 
$0$-th term of $\tau_{2m}\Omega_{\comm}$ is the sheaf cokernel of 
$d:\Omega^{2m-1}_{\comm}\to\Omega^{2m}_{\comm}$; it is not a 
quasi-coherent sheaf. There is always a map
$$
H_{dR}^{n+2m}(X)\to \Hy^n(X,\tau_{2m}\Omega_{\comm})
$$
induced by the projection $\Omega_{\comm}[2m]\fib \tau_{2m}\Omega_{\comm}$ but it is not an isomorphism
in general, not even if $X$ is affine. Hence in general
$$
\Hy^*(\Spec R,\frac{\Omega_{\ncl}}{[\Omega_{\ncl},\Omega_{\ncl}]})\ne 
H^*(\frac{\Omega_{\ncl}R}{[\Omega_{\ncl}R,\Omega_{\ncl}R]})\tag{6}
$$
 Next we consider a third type of de Rham complex; the periodic
$\fx$-complex of [CQ2]. Recall that if $R$ is any algebra then $\fx R$ is the $2$-periodic complex
with 
$$
\fx^{\op{even}}R=\Omega^0R=R\qquad \fx^{\op{odd}}R=\Omega^1R_\natural:=\frac{\Omega^1R}{[R,\Omega^1R]}
$$
and with the de Rham differential as coboundary from even to odd degree and the Hochschild boundary
from odd to even degree. We compute the cohomology of $\fx$ for formally $\nci$-smooth schemes. Such a 
gadget consists of a commutative scheme $X$ together with an inverse system of Zariski sheaves
$$
\calo_{{Y_\infty}}:\dots\fib\calo_{Y_{l}}\fib\calo_{Y_{l-1}}\fib\dots\fib\calo_{Y_{1}}
\fib\calo_{Y_{0}}=\calo_X
$$
such that each $Y_l=(X,\calo_{Y_l})$ is a formally $\ncl$-smooth scheme, that
$$
\frac{\calo_{Y_{l}}}{F_l\calo_{Y_{l}}}=\calo_{Y_{l-1}}
$$
and that the map $\calo_{Y_{l}}\fib \calo_{Y_{l-1}}$ is the natural projection. 
We compute the hypercohomology of the procomplex $\fx\calo_{Y_\infty}:=\{\fx\calo_{Y_l}\}_l$; 
we show (Corollary 6.10) that
$$
\Hy^n(X_{\pro-\zar},\fx(\calo_{Y_{\infty}}))=\prod_{2j\ge n}H_{dR}^{2j-n}X\tag{7}
$$
Moreover the decomposition above is induced by the commutator filtration. On the other hand
we prove that for periodic cyclic homology
$$
HC^{per}_* (X)=\Hy^*(X_{\pro-\zar},\fx\calo_{{Y_\infty}})\tag{8}
$$
Putting \thetag{7} and \thetag{8} together we get the well-known formula ([L, 5.1.12], [W3, 3.3])
$$
HC^{per}_n(X)=\prod_{2j\ge n}H_{dR}^{2j-n}X\tag{9}
$$
Our proof gets rid of the usual finiteness hypothesis and shows that the so-called Hodge decomposition 
\thetag{9} comes from the commutator filtration.

Each of the results for formally smooth schemes mentioned up to here is deduced from a general theorem
which holds without any smoothness hypothesis. In the absence of formal smoothness we need to replace
de Rham by infinitesimal cohomology. Recall that
if $X$ is a commutative scheme then its infinitesimal cohomology is the cohomology of the structure
sheaf on the infinitesimal site, which consists of all nilpotent thickenings $U\harrow T$ of
open subschemes $U\subset X$. One can also consider the $\ncl$-infintesimal site of any $\ncl$-scheme $Y$, consisting of all nilpotent
thickenings $U\harrow T$ of open subsets $U$ of $Y$ with $T$ an $\ncl$-scheme. 
We prove that for $Y$ formally $\ncl$-smooth 
$$
\align
H^*(Y_{\ncl-\inf},\calo)&=\Hy^*(Y_{\zar},\Omega_{\ncl})\tag{10}\\
H^*(Y_{\ncl-\inf},\frac{\calo}{[\calo,\calo]})&=\Hy^*(Y_{\zar},\frac{\Omega_{\ncl}}{[\Omega_{\ncl},
\Omega_{\ncl}]})\tag{11}
\endalign
$$
The above generalizes the theorem of Grothendieck ([Dix]) which establishes the case $l=0$ of \thetag{10} 
(compare also [Co2, Th. 3.0]). For commutative but not necessarily 
$\comm$-formally smooth $X$, 
we have (as part of theorems 6.2 and 6.5) the following generalizations of \thetag{3} and \thetag{5}
$$
\align
H^*(X_{\ncl-\inf},\calo)=&H^*(X_{\comm-\inf},\calo)\\
H^*(X_{\ncl-\inf},\frac{\calo}{[\calo,\calo]})=&
\bigoplus_{m=0}^lH^{*+2m}(X_{\comm-\inf},\frac{\Omega_{\comm}^m}{d\Omega_{\comm}^{m-1}})
\endalign
$$
In place of \thetag{7} and \thetag{8}  we obtain (as part of Theorem 6.8)
$$
\align
HC^{per}_n(X)
=&\Hy^n(X_{\nci-\inf},\fx)\\
=&\prod_{2m\ge n}H^{2m-n}(X_{\comm-\inf},\calo)\qquad (n\in\Bbb Z)
\endalign
$$
Here $\nci-\inf$ denotes the site of all nilpotent thickenings $U\harrow T$ with $T$ an $NC$-scheme of
arbitrary order. In particular we recover Feigin-Tsygan's formula ([FT,Th. 5], [W3, Th. 3.4]) 
without finiteness hypothesis and by noncommutative methods, showing that also this instance of the 
Hodge decomposition comes from the commutator filtration. 
In the commutative case there is an equivalent definition of infinitesimal cohomology which also 
generalizes to $NC$-schemes and is as follows. Assume first that $X$ admits a closed embedding 
$\iota:X\harrow Y$
into a formally smooth scheme $Y$, with ideal of definition $I\subset \calo_Y$. Then it is known
that the infinitesimal cohomology of $X$ is the same thing as the hypercohomology of the $I$-adic
completion 
$$
\Hy^*(X_{\zar},\iota^{-1}\widehat{\Omega}_{\comm}\calo_Y)=H^*(X_{\comm-\inf},\calo)\tag{12}
$$
The same is true in the case of $\ncl$-schemes, with 
$\Omega_{\ncl}$ and a formally $\ncl$-smooth scheme $Y_l$ substituted
for $\Omega_{\comm}$ and $Y$ (cf. Theorem 6.2). A similar statement holds for 
$\calo/[\calo,\calo]$ and $\Omega_{\ncl}/[\Omega_{\ncl},\Omega_{\ncl}]$, but we have to take pro-complex
cohomology rather than just complete (cf. Theorem 6.5). 
Back to the commutative theory, when $X$ cannot be embedded in a formally smooth
scheme, one can still take an open covering $\calu$ of $X$ consisting of embeddable schemes 
(affine schemes are embeddable). Then one can combine the completed de Rham complexes of each of the
local embeddings and of their intersections into a kind of \v Cech complex as done in [H1, p. 28]; 
the analogy of \thetag{12} holds for this complex. The 
same is true in the $NC$-case, and \thetag{12} as well as its version for
$\calo/[\calo,\calo]$ hold for systems of local $NC$-embeddings (defined in 3.2 below) substituted for 
single $NC$-embeddings (Theorems 6.2 and 6.5). 

We also consider the Poisson analogue of the commutator filtration, obtained by substituting Poisson
for Lie brackets in \thetag{1}. This leads one naturally to the notion of $\npl$-schemes,
their differential forms, their nilpotent thickenings and through the latter to their infinitesimal
topologies. We show that if $X$ is a commutative scheme then  
$$
\align
H^*(X_{\npl-\inf},\calo)=&H^*(X_{\ncl-\inf},\calo)\\
H^*(X_{\npl-\inf},\frac{\calo}{\{\calo,\calo\}})=&H^*(X_{\ncl-\inf},\frac{\calo}{[\calo,\calo]})
\tag{13}\\
H^*(X_{\npl-\inf},\fy)=&H^*(X_{\ncl-\inf},\fx)\\
\endalign
$$
Here we write $\calo$ for the structure sheaf of both the $\ncl$- and $\npl$-infinitesimal sites,
and $\fy$ is a Poisson adaptation of the $\fx$-complex, similar to the adaptation of the usual cyclic
complex given in [Br] and [Kas] (see 6.7-8 below).

The rest of this paper is organized as follows. In section 2 the notion of a quasi-coherent sheaf
of bimodules on an $NC$-scheme is introduced, and its elementary properties are proved. Then this is
used to establish the $NC$-analogues of several notions from elementary algebraic geometry and their
basic properties.
In section 3 the $\ncl$-infinitesimal site of a scheme is introduced. The connection between
this site and the indiscrete infinitesimal site of an algebra considered in [Co1] and [Co2] is
discussed (3.4-3.5). This section also contains a useful
lemma regarding the \v Cech-Alexander complex for infinitesimal topology (Lemma 3.3.1). Section 4
concerns $NC$-differential forms and their elementary properties. 
In section 5 the Poissson analogues
of what has been done in previous sections are discussed.  
In section 6 the main results of the paper are stated. These are packed into three theorems.
The first (6.2) computes the
$NC$- and $NP$-infinitesimal cohomologies of the structure sheaf, and the Zariski
hypercohomology of the complexes of $NC$- and $NP$-forms. The second (6.5) computes the cohomology
of the structure sheaves modulo Lie and 
Poisson brackets and compares them with the hypercohomology of the complexes of forms modulo commutators
and Poisson brackets. 
The third (6.8) computes the infinitesimal hypercohomology
of the Cuntz-Quillen complex and of its Poisson analogue. All three theorems are stated in their fullest 
generality; for the reader's convenience the particular case of each theorem concerning formally smooth 
schemes has been included as a corollary. The proofs of the main theorems are given in section 8, after
a number of lemmas and auxiliary results which are the subject of section 7. Among these auxiliary
results at least one is of independent interest (Prop. 7.9). It establishes that if $X$ is a commutative
scheme then
$$
H^p(X_{\comm-\inf},\Omega_{\comm}^q)=0\qquad \text{ for } p\ge 0, q\ge 1
$$ 
\bigskip
\head{2. Basic properties of $NC$-schemes}\endhead
\subhead{2.1. $NC$-bimodules and associated sheaves}\endsubhead
We extend the commutator filtration \thetag{1} to arbitrary $R$-bimodules $M$ by $F_{0}M=M$ 
and
$$
F_{n+1}M:=\sum_{p=0}^nF_pMF_{n+1-p}R+F_{n+1-p}RF_pM+<[F_pM,F_{n-p}R]>\tag{14}
$$
As in [Kap] we write $NC_{l}$ for the category of those 
algebras $R$ such that $F_{l+1}R=0$; in addition we put $\ncm(R)$ for those $M\in R-\op{Bimod}$ such that $F_{m+1}M=0$.
Let $\nci=\cup_{l\ge 0}\ncl$, $\nci(R)=\cup_{m\ge 0}\ncm(R)$. Note
that $R\in\ncl\Rightarrow R\in\ncl(R)$ ($0\le l\le \infty$). If $\alpha:R@>>>R'\in\nci$ is a
homomorphism and $N$ is an $R'$-bimodule then we can either take its
commutator filtration as an $R'-$ or an $R-$ bimodule (via $\alpha$). For each $0\le n$ we
have the inclusion
$$
F^R_n(N)\subset F^{{R'}}_n(N)\tag{15}
$$
It follows that 
$$
N\in NC_n(R')\Rightarrow N\in NC_n(R)\tag{16}
$$
In other words the functor $R'-\op{Bimod}@>>> R-\op{Bimod}$ induced by $\alpha$ sends
$NC_n(R')$ into $NC_n(R)$ ($0\le n\le\infty$). 
For $n<\infty$ the functor $NC_n(R')@>>>NC_n(R)$ has a left adjoint given by
$$
M\mapsto
\frac{R'\otimes_RM\otimes_RR'}{F_{n+1}(R'\otimes_RM\otimes_RR')}
$$
Note however, that the functor $\nci (R')\to \nci (R)$ does not have a left 
adjoint. 

Fix $R\in NC_l$ and $M\in NC_m(R)$; put $A:=R/F_1R$. The proof of
[Kap,2.1.5] shows that any multiplicative subset $\hat{\Gamma}\subset R$
satisfies both the right and the left \O re conditions. Similarly the
proof
of [Kap, 2.1.7] shows that localization commutes with each of the terms of
the filtration \thetag{1}. In particular
$F_1(R[\hat{\Gamma}^{-1}])=(F_1R)[\hat{\Gamma}^{-1}]\subset R[\hat{\Gamma}^{-1}]$ is a nilpotent
ideal, whence an element of $R[\hat{\Gamma}^{-1}]$ is invertible if it is so in
$R[\hat{\Gamma}^{-1}]/(F_1R)[\hat{\Gamma}^{-1}]$ which --by exactness of \O re
localization-- is the same thing as $\Gamma^{-1}A$, the commutative localization at the
image $\Gamma\subset A$ of $\hat{\Gamma}$. We have just shown that $R[\hat{\Gamma}^{-1}]$
depends only on $\Gamma$. We shall therefore write $R[\Gamma^{-1}]$ to mean
$R[\hat{\Gamma}^{-1}]$. The identity
$$
s^{m+1}x=(\sum_{i=0}^ms^{m-i}ad(sa)^i(x))s\qquad (s\in R, x\in M)
$$
is proved in the same manner as [Kap 2.1.5.1]. It shows that right
multiplication by $s$ is surjective if left multiplication is. Similarly,
by the same argument as in {\it loc. cit.}, $xs=0$ implies $s^{m+1}x=0$,
whence also injectivity of left and right multiplication are equivalent. 
It follows that there is a canonical isomorphism
$$
M[\Gamma^{-1}]:=R[\Gamma^{-1}]\otimes_R M\cong M\otimes_R R[\Gamma^{-1}]
$$
The same proof as in [Kap 2.1.6] applies with $M$ substituted for $R$ and proves
that 
$$
F_n(M[\Gamma^{-1}])=(F_nM)[\Gamma^{-1}]\qquad (n\ge 0)\tag{17}
$$
Here the $F_n$ on the left is taken in the sense of $R[\Gamma^{-1}]$-bimodules,
while that on the right is taken in $R-\op{Bimod}$. In particular 
$M[\Gamma^{-1}]\in NC_m(R[\Gamma^{-1}])$. Next we show that $M$ and its
localizations define a sheaf $\tilde{M}$ on $\Spec A$. We need some
notations. If $f\in A$, $\fp\in \Spec A$ and $x\in M$ we put 
$$
D(f)=\{\fq\in \Spec A: f\notin \fq\}
$$
and write $M_f$ and $M_{\fp}$ for the localizations at $\{f^n:n\ge 0\}$ and
at $A\backslash \fp$ respectively and $x_{\fp}$ for the image of $x$ in $M_{\fp}$.
If $U\subset \Spec A$ is open, we put 
$$
\tilde{M}(U)\subset\prod_{\fp\in U}M_{\fp}
$$
for the subset of all those elements $\sigma$ which satisfy
$$
(\forall \fp\in U)(\exists U\supset D(f)\owns \fp, x\in
M_f)(\forall\fq\in D(f))\quad \sigma_{\fq}=x_{\fq}\tag{18}
$$
It turns out that 
$$
\tilde{M}(D(f))=M_f\tag{19}
$$
Seeing this amounts to showing that the sheaf condition holds for affine
coverings of affine open subsets, and one reduces immediately to the case
when the affine open is all of $\Spec A$. For $M=R$, this is [Kap,2.2.1(b)]. 
In view of \thetag{17} the same proof as in {\it loc. cit.}
works for arbitrary $M$. We put 
$$
X=\Spec R:=(\Spec A, \tilde{R})
$$
and $\calo_X:=\tilde{R}$. Note that $\Spec R$ is a locally ringed space (in the
obvious noncommutative sense) and that $\tilde{M}$ is an $\calo_X$-bimodule.
Note further the adjoint property
$$
\hom_{\calo_X-\op{Bimod}}(\tilde{M},\calg)=\hom_{R-\op{Bimod}}(M,\calg (X))\tag{20}
$$
for $\calg \in\calo_X-\op{Bimod}$. The $\calo_X$-bimodule $\tilde{M}$ is an example
of the general notion of $NC_m$-bimodule over $\calo_X$, which is defined
as follows. If $\Cal G$ is an $\calo_X$-bimodule, we write $F_n\Cal G$ for the 
sheafification of the presheaf $U\mapsto F_n\Cal G(U)$ ($n\ge 0$), with 
$F_n(\Cal G(U))$ taken in the sense of $\calo_X(U)$-bimodules. Note the
inclusion
$$
F_n(\Cal G(U))\subset(F_n\Cal G)(U)\tag{21}
$$
We say that $\Cal G$ is an $NC_m$-{\it bimodule over} $\calo_X$ and write
$\Cal G\in NC_m(\calo_X)$ if $F_{m+1}\Cal G=0$. For $\Cal G=\widetilde{F_nM}$ the
inclusion \thetag{21} together with the adjunction \thetag{20} give
a sheaf map
$$
\widetilde{F_nM}@>\cong>>F_n\tilde{M}\tag{22}
$$
which is an isomorphism by \thetag{17} and \thetag{19}. In particular $\tilde{M}\in
NC_m(\calo_X)$. We remark that the functor $NC_m(R)@>>>NC_m(\calo_X)$ which sends
$M\mapsto \tilde{M}$ is exact, because \O re localization is. Suppose now
a homomorphism $\alpha:R@>>>R'\in \nci$ is given. Then $\alpha$ descends
to a homomorphism $A@>>>A'=R'/F_1R'$, which in turn induces a continuous
map $\hat{\alpha}:\Spec A'@>>>\Spec A$. The map $\hat{\alpha}$ together with
the induced homomorphisms
$$
\tilde{R}(D(f))=R_f@>>>R'_{\alpha
f}=\hat{\alpha}_*(\tilde{R'})(D(f))\qquad 
(f\in A)
$$
give rise to a map of locally ringed spaces $X':=\Spec R'@>>>X$. If 
$\calg \in NC_m(\calo_{X'})$ then $\hat{\alpha}_*\calg\in NC_m(\calo_X)$, by
\thetag{16}. If furthermore $\calg=\tilde{N}$ for some $N\in NC_m(R')$
then 
$$
\hat{\alpha}_*\calg=_\alpha\negthickspace\tilde{N}_\alpha
$$
where the subscript indicates that $R$ acts through $\alpha$. A left
adjoint of the functor $\hat{\alpha}_*:NC_m(\calo_{X'})@>>>NC_m(\calo_{X})$ is given
by 
$$
\hat{\alpha}^*\calg=\frac{O_{X'}\otimes_{\hat{\alpha}^{-1}\calo_X}\otimes
\calg \otimes_{\hat{\alpha}^{-1}\calo_X}O_{X'}}
{F_{m+1}(O_{X'}\otimes_{\hat{\alpha}^{-1}\calo_X}\otimes
\calg \otimes_{\hat{\alpha}^{-1}\calo_X}O_{X'})}\tag{23}
$$ 
\bigskip
\subhead{2.2. Quasi-coherent sheaves}\endsubhead
Let $X=(X,\calo_X)@>>>\Spec k$ be a (not necessarily commutative) locally ringed space over $\Spec k$ and $0\le l\le \infty$. 
We say that $X$ is an {\it affine $\ncl$-scheme} if it is isomorphic --as a locally ringed space over $\Spec k$-- to the 
spectrum of some $R\in\ncl$, and in general that it is an $\ncl$-{\it scheme} if every point $p\in X$ has an open neighborhood
$U$ such that $(U,\calo_{X|U})$ is an affine $\ncl$-scheme. We write $\ncl-Sch$ for the category of $\ncl$-schemes and
morphisms of locally ringed spaces and put $\nci-Sch=\cup_{l\ge 0}\ncl-Sch$. Note that $NC_0-Sch$ is the usual category
of commutative schemes over $\Spec k$. Like in the commutative case, the global sections functor is right adjoint
to $\Spec$; we have
$$
\Hom_{\ncl-Sch}(X,\Spec R)=\Hom_{\ncl}(R,\calo_X(X))
$$
This is proved in two steps, first for $X$ affine and then in general; the arguments of the proofs of [EGA, 1.7.3]
and [EGA, 2.2.4] apply verbatim to the $NC$-case.

Fix $X\in\ncl$; an $\ncm$-bimodule over $\calo_X$ as defined in 2.1 above is called {\it quasi-coherent} if every point
$p\in X$ has an open affine neighborhood $U$ such that the natural map
$$
\widetilde{M(U)} @>\sim>>M_{|U}\tag{24}
$$
is an isomorphism. Put $QCoh_m(X)$ for the category of $\ncm$-quasi-coherent bimodules. Recall [H2, Prop. 5.4] that
for $l=0$ the definition we have just given is equivalent to the condition that \thetag{24} be an isomorphism for
{\it every} affine open subset $U$. The same is true for arbitrary $m$ and $l$. To see this note that 
in 2.1 we have
already proved the analogues of Prop. 5.1 and 5.2 and of Ex. 5.3 of {\it loc. cit.} One checks, using 
these results,
together with elementary properties of \O re localization, that the proof of Lemma 5.3 in {\it loc. cit.} goes through
for arbitrary $m$ and $l$. The $\ncm$-analogue of Prop. 5.4 of {\it loc. cit.} is then immediate. As an application
of all this as well as of \thetag{22} and of the exactness of the functor $\tilde{}$  we get that if $M\in QCoh(X)$
and $U$ is affine then
$$
F_nM_{|U}\cong\widetilde{F_nM(U)}\text{ and } \frac{M}{F_nM}{|_U}=\widetilde{\frac{M(U)}{F_nM(U)}}\qquad\text{($U$ affine)}
\tag{25}
$$
In particular for each $0\le n\le l<\infty$, the locally ringed space
$$
X^{[n]}:=(X,\calo_X),\qquad \calo_{X^{[n]}}:=\frac{\calo_X}{F_{n+1}\calo_X}
$$
is an $NC_n$-scheme. We have a canonical identification
$$
QCoh_m(X)=QCoh_m(X^{[n]})\qquad (0\le m\le n)
$$
If $\cals$ is any abelian sheaf on the commutative scheme $X^{[0]}$, we put
$$
H^*(X_{\zar},\cals):=H^*(X^{[0]}_{\zar},\cals)\tag{26}
$$
where the subscript indicates that cohomology is taken with respect to the Zariski topology.
If $M$ is $\ncm$-quasi-coherent then the commutator filtration induces a cohomology spectral sequence 
$$
E_2^{p,q}=H^p(X_{\zar}^{[0]},\frac{F_qM}{F_{q+1}M})\Rightarrow H^{p+q}(X_{\zar},M)
$$
We remark that --by \thetag{25}-- the schemes $F_nM/F_{n+1}M$ are quasi-coherent $(n\ge 0)$. Thus for example
if $X$ happens to be affine then
$$
H^n(X_{\zar},M)=0\qquad (n>0)\tag{27}
$$
As an application of \thetag{27} one obtains that the subcategory 
$$
QCoh_m(X)\subset\ncm(\calo_X)
$$ 
is closed under extensions; indeed the proof of [H2, Prop. 5.7] applies. Now let $f:X@>>>Y\in\nci-Sch$ be a homomorphism.
The inclusion \thetag{15} implies that the functor $f_*:\calo_X-\op{Bimod}@>>>\calo_Y-\op{Bimod}$ sends 
$\ncm(\calo_X)$ into $\ncm(\calo_Y)$. Formula \thetag{23} defines a left adjoint functor $f_m^*$
of the induced functor $f_*^m:\ncm(\calo_X)@>>>\ncm(\calo_Y)$. It is clear from the affine case 
\thetag{2.1} that $f_m^*$ always sends $QCoh_m(Y)$ into $QCoh_m(X)$. The proof of [H2, 5.8 c)] shows 
that if $X^{[0]}$ is noetherian, then also $f_*^m$ preserves quasi-coherence.

We say that the morphism $f$ is a {\it closed} or an {\it open immersion} if it is so in the sense of locally ringed
spaces. The argument of the proof of [EGA, 4.2.2-b)] shows that if $f$ is a closed immersion and $Y$ is an affine 
$\ncl$-scheme then also $X\in\ncl$ and is affine. One shows using this that for any closed immersion $f$ the functor
$$
f_*^m:\ncm(\calo_X)@>>>\ncm(\calo_Y)\tag{28}
$$
preserves quasi-coherence. As an application one obtains a one-to-one correspondence between equivalence classes of closed
immersions $X\harrow Y$ and quasi-coherent two sided ideals of $\calo_Y$. 
\bigskip
\proclaim{Lemma 2.2.1} Let $X\in\nci$, $M@>g>>N\in QCoh_l(X)$, $\bar{g}:\bar{M}:=M/F_1M@>>>\bar{N}$ the induced map. Then
\item{i)} $g$ is surjective $\iff$ $\bar{g}$ is.
\item{ii)} Assume $M=N$. If $\bar{g}=id_{\bar{M}}$ then $g_{|F_lM}=id_{F_lM}$.
\endproclaim

\demo{Proof} Part $\Rightarrow$ of \thetag{i} is trivial. To prove the converse we may assume $X$ affine. Furthermore
by \thetag{27} it suffices to show that if $R\in\nci$ and $h:P@>>>Q\in\ncl(R)$ is such that $\bar{h}$ is surjective
then so is $h$. To prove this it suffices to show $h_n:G_nP:=F_nP/F_{n+1}P@>>>G_nQ$ is surjective for all $n\ge 1$. 
Every element of $G_nQ$ is represented by a sum of elements of the form $a\cdot x\cdot b$ where 
$$
x=ad(r_1)\circ\dots\circ ad(r_j)(q)\qquad q\in Q, r_1,\dots,r_j\in R, a\in F_iR, b\in F_kR, i+j+k=n
$$
If $h(p)\equiv q\mod F_1Q$, then
$$
h(F_nP)\owns h(a\cdot  ad(r_1)\circ\dots\circ ad(r_j)(p)\cdot b)\equiv a\cdot x\cdot b\mod F_{n+1}Q
$$
This proves \thetag{i}. To prove \thetag{ii} assume 
that $P=Q$ and that $\bar{h}$ is the identity. Then for $n,a,x$ and $b$ as above, 
$$
\align
h(a\cdot x\cdot b)=&a\cdot ad(r_1)\circ\dots\circ ad(r_j)(h(q))\cdot b\\
                  \equiv & a\cdot x\cdot b\mod F_{n+1}P\\
\endalign
$$
For $n=l$, $F_{n+1}P=0$, so $\equiv$ can be replaced by $=$.\qed
\enddemo
\bigskip
\proclaim{Corollary 2.2.2} Let $f:X@>>>Y\in\nci-Sch$. Then $f$ is a closed immersion $\iff$ $f^{[0]}$ is.\endproclaim
\bigskip
\remark{Remark 2.2.3} For open instead of closed immersions we still have 
$$
f:X@>>>Y\text { open immersion }\Rightarrow f^{[n]}:X^{[n]}@>>>Y^{[n]}\text{ open immersion } 
$$
for each $n\ge 0$.
\endremark 
\bigskip
\subhead{2.3. Products of $NC$-schemes; separated schemes}\endsubhead
For $l<\infty$ the categorical product of two affine $\ncl$-schemes as objects of $\ncl-Sch$ is given by 
$$
\Spec R\times_l\Spec R'=\Spec \frac{R*R'}{F_{l+1}(R*R')}\tag{29}
$$
where $*$ is the coproduct in the category $\Ass$ of associative algebras. The product of not necessarily affine 
$X,X'\in\ncl-Sch$  --denoted $X\times_l X'$-- is constructed by glueing together products of affine ones, just as
in the commutative case. Note that products do not exist in $\nci$.
We say that an $NC_l$ scheme is {\it separated} 
--over $\Spec k$-- if the diagonal map $\delta_l:X@>>>X\times_l X$ is a closed immersion.

\bigskip
\proclaim{Lemma 2.3.1}Let $X,Y\in\ncl$, $\infty>l\ge m\ge 0$. Then 

\item{i)} $(X\times_lY)^{[m]}=X^{[m]}\times_m Y^{[m]}$.

\item{ii)} $X$ is separated $\iff$ $X^{[0]}$ is separated.
\endproclaim

\demo{Proof} To prove part \thetag{i}. The projections $X\times_l Y\rightarrow X,Y$ induce a map 
$f_m:(X\times_l Y)^{[m]}@>>>X^{[m]}\times_m Y^{[m]}\in\ncm-Sch$. To show $f_n$ is an isomorphism we may assume
$X,Y$ are affine, in which case the lemma is immediate from \thetag{29}. Part \thetag{ii} is immediate from 
\thetag{i} and corollary 2.2.2.\qed
\enddemo
\bigskip

\subhead{2.4. Thickenings}\endsubhead
In this paper by a {\it thickening} of an $\nci$-scheme $X$ we understand  
a closed immersion $\tau:X@>>>T\in\nci-Sch$ such that 
$J_\tau:=\ker(\calo_T\fib\tau_*\calo_X)$ is a nilpotent ideal. 
If both $X,T\in\ncl-Sch$ then we say that $\tau$ is an $\ncl$-thickening. 
For example if $X\in\nci$ then for each $0\le m\le l$ the inclusion
$$
X^{[m]}\harrow X^{[l]}\tag{30}
$$
is an $\ncl$-thickening. We remark that all $NC$-thickenings considered in [Kap] are either of the form 
\thetag{30} or colimits
of such. However the definition given here is more general, as it includes for example all thickenings of commutative
schemes in the commutative sense ([Dix, 4.1]); indeed these are precisely the 
$NC_0$-thickenings. In fact we
have
\bigskip
\proclaim{Lemma 2.4.1} $\tau:X\to T$ is an $\ncl$-thickening $\iff$ $\tau^{[0]}:X^{[0]}\to T^{[0]}$ is an $NC_0$-thickening.
\endproclaim
\demo{Proof} Immediate from 2.2.2.\qed
\enddemo
\bigskip
\proclaim{Lemma 2.4.2}Let $\tau:X\harrow T\in\nci-Sch$ be a thickening. Then $X$ is affine $\iff$ $T$ is.
\endproclaim

\demo{Proof} If $T$ is affine then $X$ must be affine since it is closed (cf. the discussion just before
\thetag{28}). To
prove the converse, we may assume $J_\tau^2=0$. Because $\tau$ is a closed immersion, $\tau_*\calo_X$ is a quasi-coherent
$\calo-\nci$-bimodule (cf. \thetag{28}). Hence $J_\tau$ is quasi-coherent, and an object of $\nci(\calo_X)$ 
since $J_\tau^2=0$
and $\tau$ is a homeomorphism. Because $X=\Spec R$ is affine, there is an $R$-bimodule $M$ such that $J_\tau=\tilde{M}$.
Put $E=\calo_T(T)$. Taking global sections in the exact sequence
$$
0 @>>>J_\tau@>>>\calo_T@>>>\tau_*\calo_X@>>>0
$$
and using \thetag{27} we get an exact sequence of $E$-modules 
$$
0@>>>M@>>>E@>>>R@>>>0
$$
Applying the functor $\tilde{}$ to the latter sequence we get
$$
0@>>>J_\tau@>>>\tilde{E}@>>>\tau_*\calo_X@>>>0
$$
It follows that the canonical adjunction map $\tilde{E}@>>>\calo_T$ is an isomorphism.\qed
\enddemo
\bigskip
\proclaim{Corollary 2.4.3} Let $X$ be an $\nci$-scheme. Then
\item{i)} $X$ is affine $\iff$ $X^{[0]}$ is.
\item{ii)} If $U,V\subset X$ are open affine subschemes and $X$ is separated then $U\cap V$ is affine.
\endproclaim
\bigskip
\subhead{2.5. Pro-sheaves}\endsubhead If $C$ is any category, we write $\pro-C$ for the category of countably indexed
pro-objects in $C$ (cf. [CQ3], [AM]). For example if $A\in C$ is an object and 
$$
\dots\subset\Cal F_{n+1}\subset\Cal F_n\subset\dots \subset \Cal F_0= A
$$
is a filtration by subobjects, then each of the following inverse systems is an object of $\pro-C$
$$
\align
\Cal F_\infty:\dots&\harrow\Cal F_{n+1}\harrow\Cal F_n\harrow\dots \harrow \Cal F_0 \\
A: \dots & @>id>> A@>id>>A @>id>>\dots @>id>> A  \\
A/\Cal F_\infty : \dots & \to \frac{A}{\Cal F_{n+1}}\to \frac{A}{\Cal F_n}\to\dots\to\frac{ A}{\Cal F_1}
\endalign
$$
We are of course assuming that the quotients above exist in $C$. Recall $\pro-C$ is abelian if $C$ is, and
 --by [CQ3, Prop. 1.1]-- has sufficiently many
injectives if $C$ does. In particular if $X$ is an $NC$-scheme, then the category 
$\pro-ShAb(X_{\zar})$ of pro-sheaves of abelian groups has sufficiently
many injectives, and thus the right derived functors of the total global section functor 
$$
{\hat{H}}^0
:\pro-ShAb(X_{\zar})\owns\cals =\{\cals_i\}_{i\in I}\mapsto
 \lim_{i\in I}H^0(X,\cals_i)\in Ab
$$
are defined. We write
$$
H^*(X_{\pro-\zar},\cals):=R^*({\hat{H}}^0)\cals
$$
There is a cohomology spectral sequence
$$
E_2^{p,q}=H^p(X_{\zar},{\lim}^q\cals)\Rightarrow H^{p+q}(X_{\pro-\zar},\cals)\tag{31}
$$
For example if $\Cal M$ is an inverse system of quasi-coherent sheaves with surjective maps 
$$
\dots\fib\Cal M_2\fib\Cal M_1\qquad (\Cal M_n\in QCoh_\infty(X))
$$
then by \thetag{27} and [H1 Ch. 1 \S 4], the derived functors of $\lim$ vanish and we get
$$
H^*(X_{\pro-\zar},\Cal M)\cong H^*(X_{\zar},\lim \Cal M)\tag{32}
$$
An important application of pro-sheaves is to fix the problem that usual sheaf cohomology does not
commute with infinite products of abelian sheaves; that is
$$
H^*(X_{\zar},\prod_{n=1}^\infty\cals_n)\ne\prod_{n=1}^\infty H^*(X_{\zar},\cals_n)
$$
However for the pro-sheaf
$$
\prod \cals:\dots\fib\bigoplus_{n=1}^3\cals_n\fib\bigoplus_{n=1}^2\cals_n\fib\cals_1
$$
we have
$$
H^*(X_{\pro-\zar},\prod\cals)=\prod_{n=0}^\infty H^*(X_{\zar},\cals_n)\tag{33}
$$
Hypercohomology of pro-sheaves --in the Cartan-Eilenberg sense [W2, App]-- is defined in the obvious 
way, and the obvious 
generalizations of \thetag{32} and \thetag{33} hold for hypercohomology.
\bigskip
\head{3. Infinitesimal topologies}\endhead
\bigskip
\subhead{3.1. The infinitesimal topologies of an $\op{NC}$-scheme}\endsubhead
 Let $0\le l\le\infty$, $X$ an $\ncl$-scheme. The {\it $\ncl$-infinitesimal site}
on $X$ is the Grothendieck topology $X_{\ncl-\inf}$ defined as follows. The underlying category $Cat(X_{\ncl-\inf})$ 
has as objects the $\ncl$-thickenings $U\harrow T$. We write $(U,T)$ or even $T$ to mean $U\harrow T$. A map
$(U,T)@>>>(U',T')$ in $Cat(X_{\ncl-\inf})$ exists only if $U\subset U'$ in which case it is a morphism of 
$\ncl$-schemes $T@>>>T'$ such that the obvious diagram commutes. A covering of an object $T$ is a family $\{T_i\to T\}$
of morphisms such that each $T_i\to T$ is an open immersion and $\cup T_i=T$. A sheaf $\cals$ on
 $X_{\ncl-\inf}$ is the same thing as a compatible collection of Zariski sheaves 
$\{\cals_T\in Sh(T_{\zar}):T\in X_{\ncl-\inf}\}$ (cf. [BO \S 5], [Dix, 4.1]). For example the 
{\it infinitesimal structure sheaf} $\calo$ is defined by the collection $\{\calo_T\}_T$ of 
the structure sheaves of $T\in X_{\ncl-\inf}$. 
We remark that a sequence of $\ncl$-infinitesimal sheaves 
$$
0\to \cals'\to\cals\to \cals''\to 0
$$
is exact $\iff$ the sequence of Zariski sheaves 
$$
0\to \cals'_T\to \cals_T\to \cals''_T\to 0\tag{34}
$$
is exact for all $T\in X_{\ncl-\inf}$. 
An important feature of infinitesimal cohomology is that it
depends only on the underlying commutative scheme. Precisely, if $X\in\ncl-Sch$ then the inclusion
$\iota:X^{[0]}\harrow X$ is an object of $X^{[0]}_{\ncl-\inf}$; thus by composition we obtain a morphism
of topologies $F:X_{\ncl-\inf}\to X^{[0]}_{\ncl-\inf}$. With the notations of [A], we put
$$
\iota_*:=F^s:ShAb(X^{[0]}_{\ncl-\inf})\to ShAb(X_{\ncl-\inf})
$$
\bigskip
\proclaim{3.1.1. Lemma} $H^*(X_{\ncl-\inf},\iota_*\cals)=H^*(X^{[0]}_{\ncl-\inf},\cals)$.
\endproclaim
\demo{Proof}One checks that the left adjoint $\iota^*$ of $\iota_*$ is exact. On the other
hand $\iota_*$ is exact by \thetag{34}. The lemma follows from the Leray spectral sequence
associated to the morphism of topoi $\iota=(\iota_*,\iota^*):Sh(X^{[0]}_{\ncl-\inf})
\to Sh(X_{\ncl-\inf})$.\qed\enddemo
\bigskip
If $X$ is a commutative scheme and $m<l\le\infty$ then the natural inclusion
$$
\hat{j}:X_{\ncm-\inf}\harrow X_{\ncl-\inf}
$$
is a morphism of topologies. An argument similar to that of the proof of 3.1.1 shows that for 
$j_*:=\hat{j}^s$ and $\cals\in ShAb(X_{\ncm-\inf})$ we have
\bigskip
\proclaim{3.1.2. Lemma} $H^*(X_{\ncl-\inf},j_*\cals)=H^*(X_{\ncm-\inf},\cals)$\qed
\endproclaim
\bigskip
\subhead{3.2. Formal $\ncl$-smoothness; systems of embeddings}\endsubhead
Let $0\le l<\infty$. An $\ncl$-scheme $X$ is {\it formally $\ncl$-smooth} ($l<\infty$) if 
it can be covered by open affine 
schemes of the form $\Spec R$ with $R$ formally $l$-smooth in the sense of [Kap] and [Co4]. Equivalently $X$ is 
formally $\ncl$-smooth if the representable sheaf $\tilde{X}$ covers the final object $*$ of the $\ncl$-infinitesimal topos,
i.e. the map $\tilde{X}\fib *$ is an epimorphism (cf. [BO, 5.28]).
An {\it $\ncl$-embedding} of $X$ is a closed immersion $\tau:X\harrow Y$ with $Y$ formally $l$-smooth. 
If $\I=\ker(\calo_Y\fib\tau_*\calo_X)$, we consider the {\it $n$-th formal neighborhood} of $X$ along $X\harrow Y$ 
$$
(Y)_n:=(X,\tau^{-1}\frac{\calo_Y}{\I^n})\in X_{\ncl-\inf}
$$
An $\nci$-embedding
is a direct system $\caly=\{X\harrow Y_{l}\harrow Y_{l+1}\harrow\dots\}$ where $X\harrow Y_l$ is an $\ncl$-embedding  and
for $m\ge l+1$ each $Y_{m-1}\harrow Y_m$ is an $\ncm$ embedding. A {\it system 
of (local) $\ncl$-embeddings} of $X$ is a 
family  $\caly=\{\tau_i:U_i\harrow Y_i:i\in I\}$ indexed by a {\it well 
ordered set} $I$ such that 
$\calu=\{U_i:i\in I\}$ is an open covering of $X$ and each $\tau_i$ is an  
$\ncl$-embedding. The utility of the order on $I$ will be clear in 3.3 below.
The definition of a system of $\nci$-embeddings is analogous. 
\bigskip
\subhead{3.3.\v Cech-Alexander complex}\endsubhead
Let $l<\infty$, $X$ be a separated $\ncl$-scheme, 
$\caly:=\{U_i\harrow Y_i:i\in I\}$ a system of
$\ncl$-embeddings, and $\Cal S$
a sheaf of abelian groups on $X_{\ncl-\inf}$. For $i_0<\dots<i_p$ ($i_j\in I$) we 
consider the following object of $X_{\ncl-\inf}$
$$
(Y_{i_0,\dots,i_p})_n:=(U_{i_0}\cap\dots\cap U_{i_p}\harrow(Y_{i_0}\times_l\dots\times_l Y_{i_p})_n)
$$
The {\it \v Cech-Alexander (pro-)complex} of $X$ relative to $\caly$ is 
the double pro-complex of Zariski sheaves
$$
\calc_\caly^{p,q}(\cals)_n=\prod_{i_0<\dots<i_p}\cals_{(Y_{i_0,\dots,i_p}^{\times_lq+1})_n}\tag{35}
$$
with as horizontal coboundary the alternating
sum of the cofaces induced by the natural inclusions and the natural projections
$$
U_{i_0,\dots,i_p}\subset U_{i_0,\dots\underset{i_j}\to\vee,\dots,i_p}\qquad Y_{i_0,\dots,i_p}\fib
Y_{i_0,\dots\underset{i_j}\to\vee,\dots,i_p}
$$
and with as vertical coboundary the alternating sum of the cofaces induced by the 
$q+1$ distinct projections $Y_{i_0,\dots i_p}^{\times_lq+1}\fib Y_{i_0,\dots,i_p}^{\times_lq}$. In other words 
$\calc_{\caly}^{*,*}(\cals)$ is a semi-cosimplcial-cosimplicial pro-sheaf, regarded as a double cochain pro-complex in the 
usual fashion.  
Next assume $\caly$ extends to a system of $\nci$-embeddings $\calz$ and let $\calg$ be an abelian
sheaf on $X_{\nci-\inf}$. By definition $\calz$ is a sequence 
$\caly_l\harrow\caly_{l+1}\harrow\dots$ of systems of embeddings and
compatible maps. Hence it gives rise to the following
double complex in $ShAb(X)^{\Bbb N\times\Bbb N}$ 
$$
\calc^{*,*}_{\calz}(\calg)_{m,n}
:=\calc^{*,*}_{{{\caly}_m}}(\calg_{|X_{\ncl-\inf}})_{n}\qquad (m\ge l, n\ge 0)\tag{36}
$$
As a pro-object, \thetag{36} is isomorphic to the inverse system
$$
\calc^{*,*}_{\calz}(\calg)_{m}:=\calc^{*,*}_{\calz}(\calg)_{m+l,m}\qquad (m\ge 0)
$$
\bigskip
\proclaim{3.3.1.Lemma} Let $0\le l\le \infty$, $X$ a separated $\ncl$-scheme, 
$\caly=\{U_i\harrow Y_i: i\in I\}$ a system of formally
$\ncl$-smooth embeddings and $\cals$ an abelian sheaf on $X_{\ncl-\inf}$. 
Assume {\rm either} of the following hypothesis holds
\item{i)} $\Cal U=\{U_i: i\in I\}$ is locally finite.
\item{ii)} ($\forall$ $T\in X_{\ncl-\inf}$) the Zariski sheaf $\cals_T$ is quasi-coherent in the sense
of 2.2 above.
\smallskip
Then with the notations of 2.5 and 3.3,
$$
H^*(X_{\ncl-\inf},\cals)=\Hy^*(X_{\pro-\zar},\calc_{\caly}(\cals))
$$
\endproclaim
\bigskip
\demo{Proof} Assume $l<\infty$.
 Consider the following objects of the infinitesimal $\ncl$-topos
$$
\widetilde{(Y_i)}_\infty:=\coli_n\widetilde{(Y_i)}_n\qquad\qquad\tilde{\caly}=\coprod_{i\in I}\widetilde{(Y_i)}_\infty
$$
where $\widetilde{(Y_i)}_n$ is the representable sheaf. Then
$$
\align
\Hy^0(X_{\pro-\inf},\calc_{\caly}(\cals))=&h^0({\hat{H}}^0\chas)\\
                                        =&\ker(\Hom(\widetilde{\caly},\cals)\to
\Hom(\coprod_{i<j}\widetilde{(Y_i)_\infty}\times\widetilde{(Y_j)_\infty},\cals))\tag{37}\\
\intertext{where the map is the difference of those induced by the two projections 
$\widetilde{(Y_i)}_\infty\times\widetilde{(Y_j)}_\infty)\to\tilde{\caly}$. Because $\tilde{\caly}\fib *$ is an effective 
epimorphism, \thetag{37} equals}
=&H^0(X_{\ncl-\inf},\cals)\\
\endalign
$$
It remains to show that $\Hy^*(X_{\pro-\zar},\calc_{\caly}(\cdot))
:\cals\mapsto \Hy^*(X_{\pro-\zar},\chas)$ is a unversal $\delta$-functor. Under
either of the hypothesis i), ii) of the lemma, the products appearing in \thetag{35} are exact.
Indeed in the case of i) this is clear, and for ii) it follows from \thetag{27} and [H1, Ch.1\S 4].
Thus we may assume $\caly$ consists of a single
embedding $X\harrow Y$. It is clear from \thetag{34} that  
$\calc_{\caly}^q(\cdot)$ is an exact functor for each $q\ge 0$,
whence $\Hy^*(X_{\pro-\zar},\calc_{\caly}(\cdot))$ is a $\delta$- functor. It remains to show that 
the functor $\calc_{\caly}^q(\cdot)$
 preserves injectives; this will follow once we show it has an exact left 
adjoint. For the product embedding $\caly^{(q)}:=\{X\harrow Y^{\times q+1}\}$ we have
$\calc_{\caly}^q(\cdot)=\calc_{\caly^{(q)}}^0(\cdot)$. Thus we may assume $q=0$. For each $T=(\tau:U\harrow T)\in X_{\ncl-\inf}$
let 
$$
n_T=\min\{m\ge 1:\ker (\calo_T@>>>\tau_*\calo_U)^m=0\}
$$
One checks that the exact functor 
$$
\pro-ShAb(X_{\zar})\owns\{\cals_n\}_n\mapsto
\left\{\bigoplus_{\Hom(T,Y_{{n_T}})}\tau_*(\cals_{{n_T}|U})\}\right\}_{T}\in ShAb(X_{\ncl-\inf})
$$
is left adjoint to $\calc_{\caly}^0(\cdot)$. This finishes the proof for  $l<\infty$; the case $l=\infty$ is proven
similarly.\qed
\enddemo
\bigskip
\subhead{3.4. The indiscrete infinitesimal topologies of an $NC$-algebra}\endsubhead
Let $0\le l\le \infty$, $A\in\ncl$. We write $\inf(\ncl/A)$ for the category of all surjective 
homomorphisms with nilpotent kernel $B\fib A$. We equip the opposite category $\inf(\ncl/A)^{op}$ with 
the {\it indiscrete topology};
this means that if $B=(B\fib A)$ then $Cov(B)$ is the set of all isomorphisms $B\iso B'\in \inf(\ncl/A)$. 
A sheaf of abelian groups on $\inf(\ncl/A)^{op}$ with this topology is the same thing as a presheaf, 
which in turn is just a covariant functor $\calg:\inf(\ncl/A)\to Ab$. Let $X=\Spec A$. With the notations
of [A] the functor 
$$
F:\inf(\ncl/A)^{op}\to X_{\ncl-\inf},\qquad B\mapsto \Spec B
$$
is a morphism of topologies, and induces a functor between the categories of sheaves of sets
$$
f_*:=F^s:Sh(X_{\ncl-\inf})\to Sh(\inf(\ncl/A)\qquad \cals\mapsto (B\mapsto H^0(\Spec B,\cals))
$$
The left adjoint of $f_*$ is the functor
$$
f^*:=F_s:Sh(\inf(\ncl/A))\to Sh(X_{\ncl-\inf}),\qquad \calg\mapsto \widetilde{\calg}
$$
where for each $T\in X_{\ncl-\inf}$, $\widetilde{\calg}_T$ is the Zariski sheaf defined by \thetag{18}.
Because $f^*$ is exact (i.e. $(f_*,f^*)$ is a morphism of topoi) we have a Leray spectral sequence 
$$
E_2^{p,q}=H^p(\inf(\ncl/A),(R^qf_*)(\cals))\Rightarrow H^{p+q}(X_{\ncl-\inf},\cals)\tag{38}
$$
Here
$$
(R^qf_*)(\cals)(B)
=H^q(\Spec B_{\zar}, \cals_{{\Spec}B})\qquad (B\in\inf(\ncl/A))\tag{39}
$$
For example if $M$ is a sheaf of $NC$-bimodules on $\inf(\ncl/A)$ then the Zariski sheaves $\cals_T$ 
($T\in X_{\ncl-\inf}$) are all quasi-coherent, whence
\thetag{39} vanishes for $q>0$, and
$$
H^*(\inf(\ncl/A),M)=
H^{*}(X_{\ncl-\inf},\tilde{M})\tag{40}
$$
If $\calg:\inf(\ncl/A)\to Ab$ is arbitrary, we still have a natural map
$$
H^*(\inf(\ncl/A),\calg)\to
H^{*}(X_{\ncl-\inf},\tilde{\calg})
$$
but this is not an isomorphism in general.
The \v Cech-Alexander pro-complex for the indiscrete topology is constructed as follows. Assume
first $l<\infty$. Given
a sheaf $\calg:\inf(\ncl/A)\to Ab$ and a presentation
$$
0\to J\to R\fib A\to 0 
$$
of $A$ as a quotient of an algebra $R\in\Ass$, we put
$$
C^p_{\ncl}(R,J,\calg)_m=\calg\left(\frac{Cyl^p(R,J)_m}{F_lCyl^p(R,J)_m}\right)\qquad (m\ge 0)
$$
where $Cyl^p(R,J)$ is the pro-algebra of [Co2]. In case  $R_l:=R/F_lR$ is formally $l$-smooth,
the procomplex $C_{\ncl}(R,J,\calg)$ computes sheaf cohomology (cf. [Co1 5.1])
$$
\align
H^*(\holi_n C_{\ncl}(R,J,\calg)_n):=&\Hy^*(\Spec k_{\pro-\zar},C_{\ncl}(R,J,\calg))\\
=&H^*(\inf(\ncl/A),\calg)\tag{41}\\
\endalign
$$
We also put
$$
C^p_{\nci}(R,J,\calg)_m=C^p_{\ncm}(R,J,\calg)_m\qquad (m\ge 0)
$$
If $R_l$ is formally $l$-smooth for all $l$ (e.g. if $R$ is quasi-free in the sense of [CQ1])
then \thetag{41} holds for $l=\infty$ as well. As a particular case of \thetag{31} (or rather
of its hypercohomology version) we obtain a spectral
sequence
$$
E_2^{p,q}=H^p({\lim}^qC_{\ncl}(R,J,\calg))
\Rightarrow H^{p+q}(\inf(\ncl/A),\calg)\qquad(0\le l\le\infty)
$$
This spectral sequence degenerates for example when $\calg$ maps surjections with nilpotent kernel to 
surjections, as $\lim^qC_{\ncl}(R,J,\calg))=0$ for $q>0$. Hence if $M$ is as in \thetag{40}
and in addition maps surjections with nilpotent kernel to surjections, then
$$
H^{*}(X_{\ncl-\inf},\tilde{M})=H^*(\lim C_{\ncl}(R,J,M))
$$
\bigskip
\remark{3.5. Remark}One can also consider the indiscrete topology on the category 
\linebreak
$\inf(\Ass/A)$
of all nilpotent extensions $B\fib A$, where $B$ runs in the category $\Ass$ of associative algebras.
It was proved in [Co2] that for $A\in \Ass$
$$
H^*(\inf(\Ass/A),\frac{\calo}{[\calo,\calo]})=HC^{per}_*A
$$
 We remark that the indiscrete infinitesimal cohomology of $A\in\ncl$ as an associative algebra does not 
agree with its cohomology as an $\ncl$-algebra. For example if $A$ is formally $\comm$-smooth then
$HC^{per}_*A$ is as in \thetag{9} while the indiscrete $\ncl$-infinitesimal cohomology
$H^*(\inf(\ncl/A),\calo/[\calo,\calo])$ is as calculated in
8.3.4 below.
\endremark
\head{4. $NC$-differential forms}\endhead
\subhead{4.1. $NC$-forms for $NC$-algebras and schemes}\endsubhead 
We make some remarks regarding the definition of $\ncl$-forms given in the introduction \thetag{2}. 
We observe that the bimodule filtration \thetag{14} is included in the $DG$-commutator filtration; we 
have
$$
F_m\Omega^pR\subset (F_m\Omega R)^p \qquad (m\ge 0)
$$ 
In particular
$$
\Omega_{\ncl}^pR\in\ncl(R)
$$
Moreover one checks that $\Omega_{\ncl}$ localizes; if $\Gamma\subset R$ is a multiplicative system then
$$
\Omega_{\ncl}(R[\Gamma^{-1}])\cong(\Omega_{\ncl}R)[\Gamma^{-1}]
$$
Thus the following $NC$-$\calo_{{\Spec}R}$-bimodules are isomorphic
$$
\widetilde{\Omega_{\ncl}R}=\Omega_{\ncl}\calo_{{\Spec} R}
$$
It follows that in general if $X\in\ncl-Sch$ then the sheaf $\Omega_{\ncl}:=\Omega_{\ncl}\calo_X$ 
is a quasi-coherent sheaf of $DG$- algebras over $\calo_X$. 
All this generalizes to the case $l=\infty$ as follows. With the notations above, put
$$
\Omega_{\nci}R=\{\Omega_{\ncl}R\}_l\in\pro-R-\op{Bimod}
$$
If $X$ is a scheme, then the $\pro$-Zariski sheaf $\Omega_{\nci}$ is defined in the obvious way.
\bigskip
\remark{4.2. Remark} With the definitions above, the functor $\Omega_{\ncl}:\ncl\to DG\ncl$ is left
adjoint to $DG\ncl\to\ncl$, $\Lambda\mapsto \Lambda^0$. Indeed for $R\in\ncl$ and $\Lambda\in DG\ncl$
$$
\multline
\Hom_{DG\ncl}(\Omega_{\ncl}R,\Lambda)=\Hom_{DG\Ass}(\Omega R,\Lambda)=\\
=\Hom_{\Ass}(R,\Lambda^0)=\Hom_{\ncl}(R,\Lambda^0)\endmultline
$$
It follows from this that $\Omega_{\ncl}R$ is formally $\ncl$-smooth in the obvious $DG$-sense $\iff$
$R$ is formally $\ncl$-smooth. On the other hand 
$$
\frac{\Omega_{\ncl}R}{F_1\Omega_{\ncl}R}=\Omega_{\comm}(\frac{R}{F_1R})
$$
for every $R\in\ncl$. Thus if $A$ is smooth commutative and $R_l$ is an $\ncl$-smooth thickening
of $A$ in the sense of [Kap, 1.6.1] then $\Omega_{\ncl}R_l$ is a $DG$-$\ncl$-smooth thickening of the
smooth $DG\comm$-algebra $\Omega_{\comm}A$.
\endremark
\bigskip
\subhead{4.3. Forms and embeddings}\endsubhead
If $R\in\Ass$ and $J\triangleleft R$ is an ideal we put
$$
\Omega_{\ncl}(R,J):=\Omega_{\ncl}(\frac{R}{F_lR+J^\infty})\in\pro-Ab
$$
for $0\le l\le\infty$.
Similarly for $l<\infty$ if $\caly=\{\tau:X\harrow Y\}$ is an $\ncl$-embedding with ideal of definition 
$J$, we put
$$
\Omega_{\ncl}^{\caly}=\Omega_{\ncl}\left(\tau^{-1}\left(\frac{\calo_Y}{J^\infty}\right)\right)\in
\pro-ShAb(X_{\zar})\tag{42}
$$
In general if $\caly=\{U_i\harrow Y_i:i\in I\}$ is a system of $\ncl$-embeddings, $\Omega_{\ncl}^{\caly}$
is the total complex of the double pro-complex whose $p$-th column is
$$
(\Omega_{\ncl}^{\caly})^{p,*}=\prod_{i_0<\dots<i_p}\iota_*\Omega_{\ncl}^{\caly^{i_0,\dots, i_p}}
$$
Here $\caly^{i_0,\dots, i_p}:=\{U_{i_0,\dots,i_p}\harrow Y_{i_0,\dots,i_p}\}$ and 
$\iota:U_{i_0,\dots,i_p}\harrow X$ is the open immersion. If $\calz=\caly_l\harrow\caly_{l+1}\harrow\dots$
is a system of $\nci$-embeddings we put
$$
(\Omega_{\nci}^{\calz})_{l,n}=(\Omega_{\ncl}^{{\caly}_l})_n
$$
\bigskip

\head{5. $NP$-algebras and schemes}\endhead

The forgetful functors going from the category of Poisson algebras to vectorspaces and to commutative
algebras have each a left adjoint, which we write respectively $\Poiss$ and $P$. We have an isomorphism
of Poisson algebras
$$
PSV=SLV=\Poiss V\tag{43}
$$
where $V$ is a vectorspace, $L:\op{Lie-Alg}\to \op{Vect}$ is left adjoint to the forgetful functor
and $S$ is the symmetric algebra. The Poisson bracket on $SLV$
is induced by the Lie bracket of $LV$. Recall from [Co4] that if $A$ is any commutative algebra then 
$PA$ carries a natural grading
$$
PA=\bigoplus_{l=0}^\infty P_l A\tag{44}
$$
such that
$$
\{P_lA,P_mA\}\subset P_{l+m+1}A,\qquad P_lA\cdot P_mA\subset P_{l+m}A
$$
In the case $A=SV$ this grading is the same as the grading
$$
SLV=\bigoplus_{m=0}^\infty S_mLV\tag{45}
$$
induced by
$$
L_0V=V,\qquad L_{m+1}V=[L_0V,L_mV]
$$
Note that this is different from the usual grading $S=\oplus_{m=0}^\infty S^m$ of the symmetric algebra.
The analogue of the commutator filtration for Poisson algebras is the {\it Poisson filtration}
defined as follows.
Let $P$ be a Poisson algebra. Put $F_0P=P$ and inductively
$$
F_{m+1}P:=\sum_{i=1}^mF_iPF_{m+1-i}P+\sum_{i=0}^m<\{F_iP,F_{m-i}P\}>\tag{46}
$$
For example
$$
F_lPA=\bigoplus_{m\ge l}P_mA
$$
The analogues of $\ncl$-algebras and schemes in the Poisson setting are called $\npl$-algebras and
schemes. 
All what has been done for $NC$-algebras and schemes translates immediately to the 
$NP$-setting. We shall not go into the details of this translation here but shall make a few remarks
about it. First of all we note that the coproduct of two Poisson algebras as such is not the same as their
coproduct as commutative algebras or tensor product. If $P=\Poiss P/KP$ and $Q=\Poiss Q/KQ$ then 
$$
P\coprod Q=\frac{\Poiss (P\oplus Q)}{<<KP,KQ>>}\tag{47}
$$
where $<<X>>$ denotes the smallest Poisson ideal containing $X$ (cf. [Co4]). Note that as $\Poiss$
has a right adjoint, it must preserve coproducts, and that \thetag{47} simply expresses this fact.
Coproducts in $\npl$
and products in $\npl-Sch$ are defined accordingly. Second of all the right definition for the $DGP$ 
of differential forms of $\Cal A\in\Poiss$ is not $\Omega_{\comm}\Cal A$, but is defined by the 
adjointness property
$$
\Hom_{DGP}(\Omega_{\Poiss}\Cal A,Q)=\Hom_{\Poiss}(\Cal A,Q^0)
$$
For example 
$$
\Omega_{\Poiss}SLV=S(L(V\oplus dV))\supsetneq S(LV\oplus dLV)=\Omega_{\comm}SLV
$$
where $dV$ and $dLV$ are intended to be meaningful notations for the graded vectorspaces
$V[-1]$ and $(LV)[-1]$. In general if $A\in\comm$, then
$$
\Omega_{\Poiss}PA=P\Omega_{\comm}A\tag{48}
$$
where the $P$ on the right hand side is the left adjoint of the forgetful functor $DGP\subset DG\comm$.
With this definition, the same considerations as to formal smoothness remarked for $NC$-algebras 
(4.2) hold in 
the $NP$-case (see [Co4,3.3]).
\bigskip
\head{6. Statement of the main theorems}\endhead
Before stating the first theorem we need some more notations.
Recall that if $\g$ is a Lie algebra and $U\g$ its universal enveloping algebra then there is an
isomorphism of vectorspaces
$$
e:S\g\iso U\g, \qquad g_1\dots g_n\mapsto\frac{1}{n!}\sum_{\sigma}sg(\sigma)g_{\sigma_1}\dots g_{\sigma n}
\tag{49}
$$
Here $\sigma$ runs among all permutations of $n$ elements. The map $e$ is called the {\it symmetrization
map}. In theorem 6.2 below we use the particular case when $\g=LV$ is as in \thetag{43}, so that $U\g=TV$,
the tensor algebra. We use a $\hat{}$ to indicate completion of a pro-sheaf; thus for example
if $\caly$ is a system of $\ncl$-embeddings of a scheme $X$ then
$$
\hat{\Omega}_{\ncl}^{\Cal Y}=\lim_n({\Omega}_{\ncl}^{\Cal Y})_n
$$
In the statement of 6.2 below we use the fact that, as follows from lemma 3.1.2, if $X$ is a commutative 
scheme and $0\le l\le \infty$ then
$$
H^*(X_{\ncl-\inf},\calo/F_1\calo)=H^*(X_{\comm-\inf},\calo)\tag{50}
$$
where the $\calo$ on the left hand side is the structure sheaf of $X_{\ncl-\inf}$ while the one
on the right hand side is that of $X_{\comm-\inf}$. The same is true with $X_{\npl-\inf}$ substituted
for $X_{\ncl-\inf}$.
\bigskip
\proclaim{6.2.Theorem} Let $X$ be a separated commutative scheme, 
$0\le l\le \infty$, $\Cal Y$, $\Cal Z$ and 
$\Cal W$  systems of local 
$\ncl$-, $\npl$- and $\comm$ embeddings of $X$. Write $\calo$ for the structure sheaf of each of the   
 infinitesimal sites on $X$. Then there is a commutative square of natural 
isomorphisms
$$
\xymatrix{H^*(X_{{\npl}-{\inf}},\Cal O)\ar[d]^{e} & 
 \Hy^*(X_{\zar},\hat{\Omega}_{\npl}^{\Cal Z})\ar[d]^{e'}\\
H^*(X_{{\ncl}-{\inf}},\Cal O)\ar[r]^{{\alpha}_l}\ar[d]^{\pi}& 
\Hy^*(X_{{\zar}},\hat{\Omega}_{\ncl}^{\Cal Y})\ar[d]^{\pi'}\\
H^*(X_{{\comm}-{\inf}},\Cal O)\ar[r]^{\alpha_0} & \Hy^*(X_{{\zar}},\hat{\Omega}^{\Cal W}_{{\comm}})}
$$
Here $e$ and $e'$ are induced by the symmetrization map \thetag{49} and each of $\pi$, $\pi'$, 
$\pi\circ e$ and $\pi'\circ e'$ by the natural projection $\calo\fib\calo/F_1\calo$ and the isomorphism
\thetag{50}.
Moreover, if we equip each of the four vertices of the top of the diagram with the filtration induced by the
corresponding commutator or Poisson filtration 
then all three edges are filtered isomorphisms.
\endproclaim
\bigskip
\proclaim{6.3. Corollary} Assume $X$ is formally smooth. Then 
 $$
\align
H^*(X_{\comm-\inf},\Cal O)=&\Hy^*(X_{\zar},\Omega_{\comm})\\
=&\Hy^*(X_{\zar}, P_{\le l}\Omega_{\comm})\tag{51}
\endalign
$$
for each $0\le l<\infty$. If moreover $X$ admits a formally $\ncl$-smooth thickening 
$X\hookrightarrow Y_l$, then
\thetag{51} equals
$$
\align
\qquad\qquad\qquad=&\Hy^*((Y_l)_{\zar},\Omega_{\ncl})=H^*((Y_l)_{\ncl-\inf},\calo)
\endalign
$$
for each  $l<\infty$.
\endproclaim
\demo{Proof} Immediate from 6.2, \thetag{26}, \thetag{48} and 3.1.1.\qed
\enddemo
\bigskip
\remark{6.4. Notation for sheaf cokernels}In theorem 6.5 and further below there appear several expressions of the form $\cals/\Cal P$ where $\cals$ is a sheaf and $\Cal P$ is a presheaf. What is meant is the sheafification of $\cals/\Cal P$. This abuse of notation saves us further decorating already involved 
symbols.
\endremark
\bigskip
\proclaim{6.5. Theorem } Let $X$, $l$, $\Cal Y$, $\Cal Z$, $\Cal W$ and $\calo$ be as in Theorem 6.2. 
Assume the underlying open coverings of $\caly$, $\calz$ and $\calw$ are locally finite. 
Then with the convention of 6.4 there is a commutative square of natural isomorphisms
$$
\xymatrix{H^*(X_{{\npl}-\inf},\frac{\calo}{\{\calo,\calo\}})\ar[d]^{\bar{e}} & 
 \Hy^*(X_{{\pro}-{\zar}},\frac{{\Omega}_{\npl}^{\Cal Z}}{\{{\Omega}_{\npl}^{\Cal Z},
{\Omega}_{\npl}^{\Cal Z}\}})
\ar[d]^{{\bar{e}}'}\\
H^*(X_{\ncl-\inf},\frac{\calo}{[\calo,\calo]})\ar[r]^{{\bar{\alpha}}_l}\ar[d]^{\gamma}& 
\Hy^*(X_{{\pro}-{\zar}},\frac{{\Omega}_{\ncl}^{\Cal Y}}{[{\Omega}_{\ncl}^{\Cal Y},
{\Omega}_{\ncl}^{\Cal Y}]})\ar[d]^{\gamma'}\\
\prod_{m=0}^lH^*(X_{{\comm}-\inf},\frac{\Omega^m}{d\Omega^{m-1}})\ar[r]^{\beta}
&\prod_{m=0}^l
 \Hy^*(X_{{\pro}-{\zar}},\tau_{m}{\Omega}^{\Cal W}_{\comm})}\tag{52} 
$$
Here $\tau_m\Omega_{\comm}$ is the complex of sheaves
$$
\tau_m\Omega_{\comm}:\quad
\frac{\Omega_{\comm}^m}{d\Omega_{\comm}^{m-1}}@>d>>\Omega_{\comm}^{m+1}@>d>>\Omega_{\comm}^{m+2}@>d>>\dots
$$
where $\Omega_{\comm}^m/d\Omega_{\comm}^{m-1}$ is in degree $0$. Each of ${\bar{\alpha}}_l$, $\bar{e}$, and  ${\bar{e}}'$ in 
\thetag{52} is induced 
by the unbarred map with the same name in Theorem 6.2; it is filtered for the respective 
commutator and Poisson filtrations.
Each of $\gamma$, $\gamma'$
is a filtered isomorphism for the commutator filtration of its source and the filtration $F'_r=\prod_{m\le r}(\cdot)$
of its target. The map $\beta$ is a product of isomorphisms
$$
\beta_m:H^*(X_{{\comm}-\inf},\frac{\Omega_{{\comm}}^m}{d\Omega_{{\comm}}^{m-1}})@>\sim>>
 \Hy^*(X_{\pro-\zar},\tau_{m}{\Omega}^{\Cal W}_{{\comm}})
$$
of which $\beta_0=\alpha_0$ is the map of Theorem 6.2.
\endproclaim
\bigskip
\proclaim{6.6. Corollary} Assume $X$ is formally $\comm$-smooth. Then
$$
\align
H^*(X_{{\comm}-\inf},\frac{\Omega^m}{d\Omega^{m-1}})&=\Hy^*(X_{\zar},\tau_m\Omega_{{\comm}})\tag{53}\\
\Hy^*(X_{\zar},\frac{P_{\le l}\Omega_{{\comm}}}{\{P_{\le l}\Omega_{{\comm}},P_{\le l}\Omega_{{\comm}}\}})
&=\prod_{m=0}^l\Hy^*(X_{\zar},\tau_m\Omega_{{\comm}})\tag{54}\\
\endalign
$$
If moreover $X$ admits a formally $\ncl$-smooth thickening $X\hookrightarrow Y_l$ then \thetag{54} equals
$$
\qquad\qquad\qquad\qquad\
 =\Hy ((Y_{l})_{\zar},\frac{\Omega_{\ncl}}{[\Omega_{\ncl},\Omega_{\ncl}]})
 =\Hy((Y_{l})_{\ncl-\inf},\frac{\calo}{[\calo,\calo]})
$$
\endproclaim
\bigskip
\remark{6.7. Notations for cyclic homology}
Let $X$ be an $\nci$-scheme. The {\it periodic cyclic homology } of $X$ is
$$
HC^{per}_*(X):=\Hy^*(X_{\pro-\zar},\Cal{CC}^{per})
$$
where $\Cal{CC}^{per}$ is the sheafification at each level of the $2$-periodic procomplex 
$$
CC^{per}=\{\left(\bigoplus_{m=0}^{n-1}\Omega^m\right)\oplus\Omega^{n}_\natural\}_n
$$
called $\theta\Omega$ in [CQ2]. In particular for the Hochschild boundary $b$,
$$
\Omega^n_\natural:=\frac{\Omega^n}{b\Omega^{n+1}}
$$
We remark that as $CC^{per}$ is $2$-periodic, 
the Cartan-Eilenberg resolution can also be taken $2$-periodic. 
Indeed the procedure for the construction of CE-resolutions described
in [W1, Proof of 5.7.2] yields periodic resolutions for periodic complexes.
At level $n=1$, $CC^{per}$ is the periodic de Rham complex
$$
\xymatrix{{\fx}:\Omega^0\ar@/^/[r]^d&\Omega^1_{\natural}\ar@/^/[l]^b}
$$
called $X$ in [CQ2]. We also consider the analogue of the latter complex for Poisson algebras,
which is defined as follows. Recall from [Br] that if 
$\Cal A$ is a Poisson algebra then there is a boundary map
$$
\multline
\delta:\Omega^*_{\comm}\Cal A\longrightarrow\Omega^{*-1}_{\comm}\Cal A,\\
\delta(p_0dp_1\land\dots\land dp_n)=\sum_{i=0}^n(-1)^{i+1}\{p_0,p_i\}dp_1\land\dots\land\widehat{dp_i}
\land\dots\land dp_n\\
 +\sum_{i<j}(-1)^{i+j}p_0d\{p_i,p_j\}\land dp_1\land\dots\land\widehat{dp_i}\land\dots\land
\widehat{dp_j}\land\dots\land dp_n
\endmultline
$$
We put
$$
\xymatrix{{\fy}:\Cal A\ar@/^/[r]^d&\Omega^1_{\delta}\Cal A\ar@/^/[l]^\delta}
$$
where
$$
\Omega^n_\delta=\frac{\Omega_{\comm}^n}{\delta\Omega_{\comm}^{n+1}}
$$
We shall abuse notation and write $\fx$ and $\fy$ for the sheafification
of $\fx$ and $\fy$ on the various topologies for schemes considered in this paper. 
If $\caly$ and $\calz$ are systems of $\ncl$- and $\npl$- embeddings then one can 
form the
procomplexes of sheaves $\fx^{\caly}$ and $\fy^{\calz}$ in the same way as was done with the
complex $\Omega^{\caly}$ in \thetag{42}.
\endremark
\bigskip
\proclaim{6.8 . Theorem} Let $n\in\Bbb Z$, $X$ , $\Cal Y$, $\Cal Z$ as in the case $l=\infty$ of theorem
6.5 above, 
and $\fx$, $\fy$ as in
6.7. Then there is a commutative diagram of natural isomorphisms
$$
\xymatrix{\prod_{2j\ge n} H^{2j-n}(X_{\comm-\inf},\calo)\ar[d]^{f_1}&\\
  \Hy^n (X_{\npi-\inf},\fy)\ar[r]^{f_2}\ar[d]^{f_3} & \Hy^n(X_{\pro-\zar},{\fy}^{\Cal Z})\ar[d]^{f_{3'}}\\
  \Hy^n (X_{\nci-\inf},\fx)\ar[r]^{f_{2'}} &\Hy^n(X_{\pro-\zar},{\fx}^{\Cal Y})\ar[d]^{f_4}\\
&HC^{per}_n(X)\\}
$$
The map $f_1$ sends the filtration  $F'_r=\prod_{m\le r}(\cdot)$ isomorphically onto the filtration induced
by the Poisson filtration \thetag{46}; both $f_2$ and $f_{2'}$ are filtered isomorphisms.
The isomorphisms $f_3$ and $f_{3'}$ are induced by the symmetrization map \thetag{49}; they map the filtration
induced by \thetag{46} isomorphically onto that induced by \thetag{1}.
\endproclaim
\bigskip
\proclaim{6.9. Corollary} (Compare [FT, Th. 5], [W3, Th. 3.4])  There is a natural isomorphism 
$$HC^{per}_n(X)\cong \prod_{2j\ge n}H^{2j-n}(X_{\comm-\inf},\calo)\qed$$
\endproclaim
\bigskip
\proclaim{6.10. Corollary} Assume $X$ is formally smooth. Then with the notations of 2.5 above
$$
\align
HC^{per}_*(X)=&\Hy^*(X_{\pro-\zar},\fy(\{\frac{P\calo}{P_{\ge n}\calo}\}_n)\tag{55}\\
\intertext{ If in addition $X$ admits a formally $\nci$-thickening $X\hookrightarrow Y_\infty$ then the group
 \thetag{55} is also isomorphic to}
=&\Hy^*(X_{\pro-\zar}, \fx \calo_{Y_{\infty}})\qed
\endalign
$$
\endproclaim
\bigskip
\head{7. Auxiliary results}\endhead
\bigskip
\proclaim{7.1. Proposition} Let $R\in\Ass$ be a quasi-free algebra, $TR$ the tensor algebra,
$JR:=\ker(TR\fib R)$, $s\in\Hom_{\Ass}(R,TR/JR^2)$ a section of the canonical projection, 
$J\triangleleft R$ an ideal and $l\ge 0$. Then 
there are maps of pro-complexes 
$$
\xymatrix{C_{\ncl}(R,J,\calo)\ar@/^/[r]^\alpha &\Omega_\ncl(R,J)\ar@/^/[l]^\beta}
$$
and homotopies $\delta:\alpha\beta@>>>1$ and $\gamma:\beta\alpha@>>>1$ all of which are natural with 
respect to 
$R$, $s$, $J$ and $l$, and interchange commutators and graded commutators. In particular $C_{\ncl}(R,J,\calo/[\calo,\calo])$ is naturally homotopic to 
$$
\Omega_{\ncl}(R,J)/[\Omega_{\ncl}(R,J),\Omega_{\ncl}(R,J)]
$$
\endproclaim
\demo{Proof} It suffices to check that the map $\alpha:=1\otimes p:\bar{C}(R,0,\calo)\mapsto\Omega R$ of the proof of 
[Co2, 2.4] which --by the proof of [Co2, 3.1]-- preserves both the $J$-adic filtration and the 
commutator subspace, preserves also the commutator filtration, and to construct a natural homotopy 
inverse for it with the same
properties. The proof that $\alpha$ preserves the commutator filtration is similar to the proof that
it preserves the commutator subspace; one just considers the action of the full symmetric group 
$\Sigma_m$ on $\bar{T}^m$ rather than only that of the cyclic group. The map $s$ of the proof of 
[Co2, 3.1] extends to a $\Sigma_m$-equivariant contracting homotopy $\theta$ of the augmented resolution 
$\bar{T}^m\fib k$. Using $\theta$ and the perturbation lemma [Co2,  2.5] one obtains a contracting
homotopy of the mapping cone of $\alpha$, $h:M^n=\bar{C}^n\oplus\Omega^nR\to M^{n-1}$ with the matricial form
 $$
h=\left(\matrix\gamma&\beta\\
               0 &0\endmatrix\right):M^p=\nyl^pR\oplus\Omega^{p-1}@>>>M^{p-1}
$$
It follows that $\beta$ is a cochain map with $\alpha\beta=1$ and that $\gamma$ is a homotopy 
$1\to\beta\alpha$. One checks, using the equivariance of $\theta$ and $s$ and the formulas of
[Co2, 2.5] that both $\gamma$ and $\beta$ preserve both the commutator subspace and 
the commutator filtration and are continuous for the $J$-adic filtration.\qed
\enddemo
\bigskip
\proclaim{Lemma 7.2} Let $A\in\comm$, $R\in\Ass$, $\pi:R\fib A$ a surjective homomorphism, 
$\pi_{l,n}:R_{l,n}:=R/F_{l+1}R+(\ker\pi)^{n+1}\fib A$ the induced map, $\Gamma\subset B:=R/F_1R$ a multiplicative
system and $\hat{\Gamma}=\pi^{-1}(\Gamma)$. Then there is a commutative diagram with horizontal isomorphism
$$
\xymatrix{R_{l,n}[\Gamma^{-1}]\ar[dr]_{\pi_{l,n}[\Gamma^{-1}]}\ar[r]^{\overset{\phi}\to{\sim}}&R[{\hat{\Gamma}}^{-1}]_{l,n}
\ar[d]^{\pi[{\hat{\Gamma}}^{-1}]_{l,n}}\\
&A[\Gamma^{-1}]}
$$
\endproclaim
\demo{Proof} Both $R\mapsto R_l[{{\Gamma}}^{-1}]$ 
and $R\mapsto (R[{\hat{\Gamma}}^{-1}])_l$ are universal (initial) among all those algebra homomorphisms going from $R$  
to an $\ncl$-algebra which invert $\Gamma$. Therefore they are isomorphic $R$-algebras. By naturality we get a commutative
diagram
$$
\xymatrix{
R_l[{{\Gamma}}^{-1}]\ar[r]^{\overset\psi\to\sim}\ar[dr]_{\pi_l[\Gamma^{-1}]}&
(R[{\hat{\Gamma}}^{-1}])_l\ar[d]^{(\pi[{\hat{\Gamma}}^{-1}])_l}\\
&\Gamma^{-1}A}
$$
where $\psi$ is the natural isomorphism of $R$-algebras just defined. Thus $\psi$ maps 
$K:=\ker(\pi_l[{{\Gamma}}^{-1}])=(\ker\pi_l)[{{\Gamma}}^{-1}]$ isomorphically to 
$K':=(\ker\pi[{\hat{\Gamma}}^{-1}])_l$.
One checks, using [Kap, 2.1.5.1] that $K^n=(\ker\pi_l)^n[{{\Gamma}}^{-1}]$. The map $\phi$ of the lemma is that
induced by $\psi$ upon passage to the quotient.\qed
\enddemo
\bigskip
\proclaim{Lemma 7.3} Let $G=G_0\oplus G_1\oplus\dots\oplus G_l$ be a graded commutative algebra. Assume $G$ is 
additionally equipped with an associative --but not necessarily commutative-- product 
$$
\Phi=\sum_{p=0}^{l}\Phi_p:G\otimes G@>>>G
$$
such that $\Phi_0$ is the original commutative product, and $\Phi_p$ is homogeneous of degree $p$ and a bidifferential
operator. Consider the associative algebra $R=(G,\Phi)$. If $\Gamma\subset G_0$ is a multiplicative system, then
$\hat{\Gamma}=\{s+g_+|s\in\Gamma, g_+\in\oplus_{n\ge 1}G_n\}\subset R$ is a multiplicative system and
$$
R[\hat{\Gamma}^{-1}]\cong (\Gamma^{-1}G,\Gamma^{-1}\Phi)
$$
\endproclaim
\demo{Proof} Note first that the product $\Gamma^{-1}\Phi$ is associative because the associator localizes
$$
\align
\Cal A(\Gamma^{-1}\Phi)=&\Gamma^{-1}\Phi(\Gamma^{-1}\Phi(,),)-\Gamma^{-1}\Phi(,\Gamma^{-1}\Phi(,))\\
=&\Gamma^{-1}\Cal A(\Phi)=0
\endalign
$$
On the other hand that $\hat{\Gamma}$ is multiplicavely closed is clear from the fact that
$\Phi(G_{\ge n}\otimes G_{\ge m})\subset G_{\ge n+m}$.
One checks that, upon localization of $G_0$-modules, the projection $G\fib G_0$ becomes a surjective algebra homomorphism
$\Gamma^{-1}\pi:R'=(\Gamma^{-1}G,\Gamma^{-1}\Phi)\fib \Gamma^{-1}G_0$. Thus an element $x\in R'$ is invertible if and only if 
$\Gamma^{-1}\pi(x)$ is invertible. It follows that the obvious homomorphism $R@>>>R'$ maps each element of $\hat{\Gamma}$
to an invertible element, whence we have a natural map $\phi:R[{\hat{\Gamma}}^{-1}]@>>>R'$. To prove that $\phi$ is an
isomorphism proceed as follows. Consider the filtration $R=R_{0}\supset R_{1}\supset\dots\supset R_{l}$, 
$R_{n}=(G_{\ge n},\Phi_{|G_{\ge n}\otimes G_{\ge n}})$. Each $R_n$ is an ideal of $R$, and by exactness of \O re localization, the associated
graded ring is
$$
\bigoplus_{n=0}^l\frac{R_n[{\hat{\Gamma}}^{-1}]}{R_{n+1}[{\hat{\Gamma}}^{-1}]}=\bigoplus_{n=0}^l\Gamma^{-1}G_n=\Gamma^{-1}G
$$
Thus $\phi$ is an isomorphism because it is so at the graded level.\qed
\enddemo
\bigskip
\proclaim{Lemma 7.4} Let $V$ be a vectorspace, $T=TV$ the tensor algebra, $T_{\le l}=T/F_{l+1}T$. Also let $P=\Poiss V$ be the
free Poisson algebra, and for the Poisson analogue of the commutator filtration, $P_{\le l}=P/F_{l+1}P$. Assume a 
multiplicative system $\Gamma\subset P_{\le 0}=T_{\le 0}=S:=SV$ is given, and let $\hat{\Gamma}\subset 
T_{\le l}$ be the inverse image of $\Gamma$ under the projection $T_{\le l}\fib S $. Then the symmetrization
map \thetag{49} induces an isomorphism 
$$
\frac{\Gamma^{-1}P_{\le l}}{\{\Gamma^{-1}P_{\le l},\Gamma^{-1}P_{\le l}\}}\cong \frac{T_{\le l}[{\hat{\Gamma}}^{-1}]}
{[T_{\le l}[{\hat{\Gamma}}^{-1}],T_{\le l}[{\hat{\Gamma}}^{-1}]]}
$$
\endproclaim
\demo{Proof} Apply lemma 7.2 with $R=T$, $A=B=S$ to obtain an isomorphism 
$T_{\le l}[\Gamma^{-1}]\cong T[{\hat{\Gamma}}^{-1}]/F_{l+1}T[{\hat{\Gamma}}^{-1}]$. Thus 
$$
\frac{T_{\le l}[{{\Gamma}}^{-1}]}{[T_{\le l}[{{\Gamma}}^{-1}],T_{\le l}[{{\Gamma}}^{-1}]]}=
\frac{T[{\hat{\Gamma}}^{-1}]}{F_{l+1}T[{\hat{\Gamma}}^{-1}]+[T[{\hat{\Gamma}}^{-1}],T[{\hat{\Gamma}}^{-1}]]}
$$
Now $[T[{\hat{\Gamma}}^{-1}],T[{\hat{\Gamma}}^{-1}]]$ is the image of the Hochschild boundary 
$b:\Omega^1T[{\hat{\Gamma}}^{-1}]@>>>T[{\hat{\Gamma}}^{-1}]$, $xdy\mapsto [x,y]$. But 
$$
\Omega^1T[{\hat{\Gamma}}^{-1}]=
T[{\hat{\Gamma}}^{-1}]dVT[{\hat{\Gamma}}^{-1}]\cong T[{\hat{\Gamma}}^{-1}]\otimes V\otimes T[{\hat{\Gamma}}^{-1}]
$$
as $T[{\hat{\Gamma}}^{-1}]$-bimodules. Hence 
$$
\frac{\Omega^1T[{\hat{\Gamma}}^{-1}]}{b\Omega^2T[{\hat{\Gamma}}^{-1}]}=
\frac{\Omega^1T[{\hat{\Gamma}}^{-1}]}{[T[{\hat{\Gamma}}^{-1}],\Omega^1T[{\hat{\Gamma}}^{-1}]]}\cong 
T[{\hat{\Gamma}}^{-1}]\otimes V
$$
and therefore
$$
[T[{\hat{\Gamma}}^{-1}],T[{\hat{\Gamma}}^{-1}]]=[T[{\hat{\Gamma}}^{-1}],V]\tag{56}
$$
On the other hand -- by [Co4 2.1 (1)] --- the map $e$ induces a vectorspace isomorphism $P_{\le l}\cong T_{\le l}$,
whence $T_{\le l}$ is identified with the algebra with underlying vectorspace $P_{\le l}$ and multiplication
$$
\Phi(x,y):=e^{-1}(exey)\tag{57}
$$
By [Co4 2.1 (2)] and [Co3 2.2], Lemma 7.3 applies to $G=P_{\le l}$ whence $T_{\le l}[{\hat{\Gamma}}^{-1}]\cong 
(\Gamma^{-1}P_{\le l},\Gamma^{-1}\Phi)$. Thus, modulo $F_{l+1}T[{\hat{\Gamma}}^{-1}]$, \thetag{56} 
gets identified
with the subspace generated by the elements of the form 
$$
\sum_{p=1}^l(\Gamma^{-1}\Phi)_p(s^{-1}x,v)-(\Gamma^{-1}\Phi)_p(v,s^{-1}x)\qquad(x\in P_{}\le l,s\in\Gamma,v\in V)\tag{58}
$$ 
By [Co3 1.1] the homogeneous part of degree one of \thetag{58} is 
$$
\{s^{-1}x,v\}\in\{\Gamma^{-1}P_{\le l},V\}
$$
I claim that each hogeneous part of degree $p\ge 2$ of \thetag{58} is zero. For $s=1$, the claim is just the fact that
$e$ commutes with the adjoint action of the Lie subalgebra $L\subset T$ generated by $V$ --[L, 3.3.5]-- and in
particular with its restriction to $V$. Recall both $\Phi_p(v,)$ and $\Phi_p(,v)$ 
are differential operators of order
$\le p$ ([Co3 2.2]). Thus the identity
$$
F(s^{-1}a)=\sum_{j=0}^p(\sum_{i=j}^p(-1)^i\binom{i}{j})s^{-(j+1)}F(s^{-1}a)\qquad (s\in \Gamma, a\in P_{\le l})
$$
holds for both $F=\Phi_p(v,), \Phi_p(,v)$. The general case of the claim follows from this observation and the case $s=1$.
We have shown that under our identifications $[T_{\le l}[{{\Gamma}}^{-1}],T_{\le l}[{{\Gamma}}^{-1}]]$ gets identified
with $\{\Gamma^{-1}P_{\le l},V\}$. It is clear that the latter coincides with 
$\{\Gamma^{-1}P_{\le l},\Gamma^{-1}P_{\le l}\}$.\qed
\enddemo
\bigskip
\proclaim{Lemma 7.5} Let $V$ be a vectorspace, $S=\Omega_{\comm}SV=S(V\oplus dV)$ the commutative $DGA$. 
Let $P_+S=\oplus_{n\ge 1}P_nS$ be the part of positive degree in the $DG$-Poisson envelope \thetag{44}.
Then there is a contracting homotopy
$$
h:P_+S@>>>P_+S
$$
which is right $S$-linear, homogeneous of degree zero for the Poisson gradation and maps $S(LV\oplus dLV)\cap P_+S$ to itself.
\endproclaim
\demo{Proof} Define a $k$-linear map $\partial:W:=V\oplus dV@>>>W$, $\partial dv=v$, $\partial v=0$. Extend $\partial$
first to $\frak{g}=LW$ as a derivation for the Lie bracket and then to all of $S\frak{g}=PS$ as a
derivation for the (skew-) commutative product. Put $\Delta=[\partial,d]$. Write $||$ for homogenous
degree with respect to \thetag{45}. Consider the grading $\omega$ of $S\frak{g}$ determined by 
$\omega(g)=|g|+1$. If $x$ is homogeneous
with respect to $\omega$, then $\Delta x=\omega(x)x$. Rescale the restriction of $\Delta$ to $S_+\g_+$
(notation as in [Co4, 1.0]) 
to obtain a $k$-linear map $\kappa:S_+\g_+@>>>S_+\g_+$ with $\kappa d+d\kappa=1$. 
$$
h:P_+S=S\g_0\otimes S_+\g_+\to P_+S,\qquad h(x\otimes y)=(-1)^{\deg x}x\otimes\kappa(y)
$$
One checks that $h$ is right $S\g_0$-linear and that 
$dh+hd=1$.\qed
\enddemo
\bigskip
\proclaim{Lemma 7.6} Let $V$ be a vectorspace,  $P=\Poiss V$, $i\ge 0$. Write $_i\Omega_{\comm}P$ for the homogeneous part of degree
$i$ with respect to \thetag{45}. Then there is a $k$-linear map 
$$
\nabla:\Omega^r_{\comm}P@>>>\Omega^{r+1}_{\comm}P\qquad{r\ge 1}
$$
such that
\item{i)} $(\nabla\delta+\delta\nabla)\omega=\omega$ if $\omega\in\Omega^r_{\comm}P$, $r\ge 2$.
\item{ii)} $\nabla(_i\Omega_{\comm}^* P)\subset_{i-1}\negthickspace\Omega^{*+1}_{\comm}P$.
\item{iii)} The restriction of $\nabla$ to $_i\Omega^r_{\comm}P$ is a differential operator of 
$SV$-modules.
\endproclaim
\demo{Proof} Put $L=LV$,
$$
C^\alpha_{r,i}:=P_i\otimes(\Lambda^rL)_{\alpha-(r+i)},\qquad C_r^\alpha:=_{\alpha-r}\negthickspace
\Omega^r_{\comm}P=\bigoplus_{i=0}^{\alpha-r}
C_{r,i}^\alpha\tag{59}
$$
We have
$$
\delta C^\alpha_{r,i}\subset \bigoplus_{j\ge i}C_{r-1,j}^\alpha
$$
whence $C^\alpha$ is a subcomplex of the complex $C=(\Omega_{\comm}P,\delta)$, and $C=\oplus_{\alpha\ge 0}C^\alpha$. The homogeneous 
component of degree $0$ of $\delta:C_r^\alpha@>>>C_{r-1}^\alpha$ is the restriction of $1\otimes\delta'$, where
$\delta':\Lambda^rL@>>>\Lambda^{r-1}L$ is the Chevalley-Eilenberg boundary
$$
\delta'(g_1\wedge\dots\wedge g_r)=\sum_{i<j}(-1)^{i+j}[g_i,g_j]\wedge g_1\wedge\dots \overset{i}\to\vee\dots\overset{j}\to\vee
\dots\wedge g_r
$$
Because $L$ is free, there is a $k$-linear map $\nabla':\Lambda^nL@>>>\Lambda^{n+1}L$ ($n\ge 1$) such that 
$\nabla'\delta'+\delta'\nabla'=1$ on $\Lambda^{\ge 2}L$. Because $\delta'$ is homogeneous of degree $+1$ for the chain complex
decomposition induced by $L=\oplus_{n\ge 0}L_n$, we may assume $\nabla'$ homogeneous of degree $-1$. Put 
$$
\nabla^{\alpha,0}:=1\otimes\nabla' :C_r^\alpha@>>>C_{r+1}^\alpha\qquad (r\ge 1)
$$
Then $\nabla^{\alpha,0}$ is homogeneous of degree $0$ for the decomposition \thetag{59}. By the perturbation lemma
([Co2 2.5]) there exists, for each $n\ge 1$, a $k$-linear map $\nabla^{\alpha,n}:C^\alpha_r@>>>C^\alpha_{r+1}$
homogeneous of degree $n$ (with respect to \thetag{59}) such that 
$$
\nabla^\alpha:=\sum_{n=0}^\infty\nabla^{\alpha,n}
$$
verifies $\nabla^\alpha\delta+\delta\nabla^\alpha=1$ on $C^\alpha_r$ ($r\ge 2$). Moreover, from the formulas of
[Co2 2.5] and the fact that each of the components of $\delta$ is a differential operator --because $\{,\}$ is 
bidifferential-- it follows that the same is true of each $\nabla^{\alpha,n}$. Therefore the map 
$\nabla=\oplus_\alpha\nabla^\alpha$ satisfies the conditions of the lemma.\qed
\enddemo
\bigskip
\proclaim{Lemma 7.7} Let $V$ be a vectorspace, $S=SV$, $P=\Poiss V$. Consider the complex
$$
\frak{N}^n=\left\{\matrix P_i&\text{ if } n=2i\\
                          \frac{_i\Omega^1_{\comm}P}{\delta_{i-1}\Omega^2_{\comm}P} &\text{ if } n=2i+1\endmatrix\right.\qquad (n\ge 0)
$$
with coboundary maps $d:\N^{2i}@>>>\N^{2i+1}$ and $\delta:\N^{2i+1}@>>>\N^{2i+2}$ ($i\ge 0$). Then $\N$ is naturally
homotopy equivalent to $\Omega_{\comm}S$ in such a way that each of the natural homotopy equivalences and homotopies involved
is continuous for the adic topology induced by any ideal $I\triangleleft S$.
\endproclaim
\demo{Proof} Put
$$
\N\Omega^n_i=\left\{\matrix \N^n &\text{ if } n\le 2i+1\\
                            \frac{_i\Omega^r_{\comm}P}{\delta_{i-1}\Omega^{r+1}_{\comm}P} &\text{ if } n=2i+r\text{ }(r\ge 2)
			    \endmatrix\right.
$$
Make $\N\Omega^n_i$ into a cochain complex with boundary map $\partial^n:\N\Omega^n_i@>>>\N\Omega^{n+1}_i$ given
by
$$
\partial^n=\left\{\matrix\delta & n \text{ odd and } \le 2i-1\\
                          d & n \text{ even or }\ge 2i
                  \endmatrix\right.
$$
Let $h$ and $\nabla$ be as in lemmas 7.5 and 7.6. Define maps
$$
\xymatrix{{\N}\Omega^n_i\ar@/^/[r]^{\alpha_n}&{\N}\Omega^n_{i+1}\ar@/^/[l]^{\beta_n}}\quad
{\N}\Omega^n_i@>\gamma_n>>{\N}\Omega^{n-1}_i\quad \text{and}\quad {\N}\Omega^n_{i+1}@>\epsilon_n>>{\N}\Omega^{n-1}_{i+1}
$$
as follows
$$
\alpha_n=\left\{\aligned 1 \qquad& n\le 2i+1\\
                  -h\delta\qquad & n\ge 2i+2\endaligned \right.\qquad \beta_n=\left\{\aligned 1\qquad&n\le 2i+1\\ 
                          -\nabla d\qquad &n=2i+2\\
                         -(\nabla d+d\nabla)\qquad & n\ge 2i+3\endaligned\right.
$$
$$
\gamma_n=\left\{\aligned 0\qquad & n\le 2i+2\\
                         \nabla h\delta\qquad & n\ge 2i+3\endaligned\right.\qquad \epsilon_n=\left\{\aligned 0\qquad
& n\le 2i+2\\
h(\delta\nabla-1)\qquad & n\ge 2i+3
\endaligned\right.
$$
One checks that $\alpha$ and $\beta$ are cochain maps as well as that the following identities hold
$$
\alpha\beta-1=\epsilon\partial+\partial\epsilon\qquad\beta\alpha-1=\gamma\partial+\partial\gamma.
$$
Thus $(\Omega_{\comm}S,d)=(\N\Omega_0,\partial)$ is naturally and adically continously homotopy equivalent to 
$\N=\coli(\N\Omega_0@>\beta>>\N\Omega_1@>\beta>>\N\Omega_2@>\beta>>\dots)$\qed
\enddemo
\bigskip
\proclaim{Lemma 7.8} Let $U$ and $V$ be vectorspaces, $\alpha,\beta\ge 0$, $\gamma=\alpha+\beta$.
Let $\Sigma_{\alpha,\beta}:=\Sigma_{\alpha}\times\Sigma_{\beta}$ act on $T^\gamma U\otimes
T^\gamma V$ as follows:
$$
\multline
(\sigma,\tau)(u_1\dots u_\gamma\otimes v_{ 1}\dots v_{\gamma})=\\
(\sg\tau) u_{\sigma 1}\dots u_{\sigma\alpha}u_{\alpha+\tau 1 }\dots u_{\alpha+\tau\beta}
\otimes v_{\sigma 1}\dots v_{\sigma\alpha}v_{\alpha+\tau 1}\dots v_{\alpha+
\tau \beta}
\endmultline
$$
Then 
$$
S^\alpha (U\otimes V)\otimes \Lambda^\beta(U\otimes V)\cong (T^\gamma U\otimes T^\gamma V)_{\Sigma_{
\alpha,\beta}}
$$
\endproclaim
\demo{Proof} Straightforward.\qed\enddemo
\bigskip
\proclaim{Proposition 7.9} Let $X$ be a separated commutative scheme. Then 
$$
H^n(X_{\comm-\inf},\Omega_{\comm}^p)=0\qquad (p\ge 1,n\ge 0)
$$
\endproclaim
\demo{Proof} By \thetag{38}, \thetag{39} and \thetag{40} it suffices to show that if $U$ is a vectorspace and $I\subset S=SU$ is an
ideal, then the normalized pro-complex $\overline{C}(S,I,\Omega^p)$ is contractible. We shall show this for the
case $I=0$; a routine verification shows that all the cochain maps and homotopies we shall define are continous for
the adic topology of any ideal $I$, proving the general case. Let $V^*$ and $W^*$ be the cosimplicial vectorspaces
of [Co2,1.2]. We have 
$$
S^{\otimes m+1}=S(U \oplus U\otimes V^m)=S(W^m)\tag{60}
$$
Hence by the lemma above 
$$
\Omega^p_{\comm}{S^{\otimes m+1}}=\bigoplus_{\alpha\ge 0}\bigoplus_{q+\beta=p}\Omega^q_{\comm}SU
\otimes(T^{\alpha+\beta}U\otimes
T^{\alpha+\beta}V^m)_{\Sigma_{\alpha,\beta}}
$$
Procompletion with respect to the ideal $<U\otimes V^m>\subset S^{\otimes m+1}$ gives the pro-space
$$
C^m_n:=C^m(S,0,\Omega^p)_n=\bigoplus_{0\le\alpha\le n}\bigoplus_{q+\beta=p}\Omega^q_{\comm}SU\otimes
(T^{\alpha+\beta}U\otimes T^{\alpha+\beta}V^m)_{\Sigma_{\alpha,\beta}} 
$$
where ($m\ge 1$). Recall from the proof of [Co2, 2.4] that for the normalized complex $\overline{T^rV}^*$ we have
$$
\overline{T^rV}^m=0\quad\text{ for }r>m\text{ and } \overline{T^mV}^m=k[\Sigma_m]
$$
Thus $\overline{C}^m$ is the constant pro-vectorspace
$$
\overline{C}^m=\bigoplus_{0\le r\le m}\bigoplus_{p-r\le q\le p}\Omega^q_{\comm}SU\otimes 
(T^rU\otimes\overline{T^rV}^m)_{\Sigma_{r-p+q,p-q}} 
$$
Recall from the proof of Proposition 7.1 that there is a $\Sigma_m$-equivariant homotopy equivalence 
$p:\overline{T^mV}@>>>k[-m]$. It follows that $1\otimes p$ passes to the quotient modulo the action of the symmetric
group, giving a homotopy equivalence between $\overline{C}$ and a complex having
$$
D^m:=\bigoplus_{p-m\le q\le p}\Omega^{q}_{\comm}SU\otimes\Lambda^{m-p+q}U\otimes S^{p-q}U
$$
in degree $m$. This vectorspace can be interpreted as a piece of the DG-module of $m$-differential forms
of the DGA $\Omega_{\comm}SU$. Namely
$$
D^m=_p\negthickspace\Omega^m_{DG-\comm}(\Omega_{\comm}SU)
$$
Here the subindex $p$ denotes weight with respect to the  
grading of $\Omega_{\comm}SU$ determined by $\deg (u)=0$, $\deg( du)=1$. One checks further that the coboundary map is the restriction of $d'$, the de Rham differential for forms on $\Omega_{SU}$. We have 
$$
\align
\Omega_{DG-\comm}\Omega_{\comm}SU=&S(U\oplus dU\oplus d'(U\oplus dU))\\
 = S(U\oplus d'U)\otimes S(dU\oplus d'dU)\cong &\Omega_{DG-\comm}SU\otimes 
\Omega_{DG-\comm}S(dU)
\endalign
$$
It is clear that $\Omega_{DG-\comm}S(dU)$ is contractible by means of a weight preserving contracting 
homotopy $h$. Thus $1\otimes h$ is a contracting homotopy
for $D$. This concludes the proof.\qed
\enddemo
\bigskip
\proclaim{Lemma 7.10} Let $n,m\ge 0$. Then 
\item{i)} Let $\Omega^1_\natural$ be the ${\nci-\inf}$ sheaf $\coker(b:\Omega^2@>>>\Omega^1)$.
Then 
$$
H^n(X_{\nci-\inf},\Omega^1_\natural)=0
$$
\item{ii)} Let $(_m\Omega^1_{\comm}P)_\delta$ be the $\comm-\inf$ sheaf 
$\coker(\delta:_{m-1}\Omega^2_{\comm}P@>>>_m\Omega^1_{\comm}P)$ where the subscript on the left hand corner indicates
degree with respect to the grading \thetag{44}. Then
$$
H^n(X_{\comm-\inf},(_m\Omega^1_{\comm}P)_{\delta})=0
$$
\endproclaim
\demo{Proof} i)It suffices to show that if $\comm\owns A=R/I$ with $R$ quasi-free then for 
the presheaf cokernel $\bar{\Omega}^*=\Omega^*/b\Omega^{*+1}$, 
the complex $C(R,I,\bar{\Omega}^1)$ is naturally contractible. A similar argument as that 
given in the proof of 7.1 above shows that the homotopy equivalence of the proof of 
[Co2, Lemma 5.6] 
$$
1\otimes p:C(R,I,\bar{\Omega}^1)
@>>>\frac{\Omega^{*+1}R\oplus\Omega^*R}{N^*+{\calg}^{1,\infty}}=
\frac{\bar{\Omega}^{*+1} R\oplus \bar{\Omega}^* R}{\calg^{1,\infty}}
$$
preserves the commutator filtration.
\item{ii)} By the proof of lemma 7.4, for $R=TV$, $S=SV$, and the presheaf cokernel 
$\bar{\Omega}^*_{\comm}:= \Omega^*_{\comm}/\delta\Omega^{*+1}_{\comm}$ the symmetrization map
induces an isomorphism of pro-complexes  
$$
\prod_{m=0}^\infty C(S,I,_m\negthickspace\bar{\Omega}^1_{\comm}P)\cong C(R,I,\bar{\Omega}^1)\qed
$$
\enddemo
\bigskip
\proclaim{Lemma 7.11} Let $V$ be a vectorspace, $S=SV$, $T=TV$, $L=LV$, $P=\Poiss V$, $I\subset S$ an ideal, 
$A=S/I$ and $J\subset T$ the inverse image of
$I$ under the projection $T\fib T/F_1T=S$. Then the map
$$
\multline
\eta:\Omega_{\comm} P@>>>\Omega T,\\
\eta(a\otimes dg_1\land\dots dg_n)=e(a)\otimes\sum_{\sigma\in \Sigma_n}\sg(\sigma) dg_{\sigma 1}\dots dg_{\sigma n}\qquad (g_i\in L)
\endmultline
$$
induces a homotopy equivalence of pro-complexes
$$
(\Omega_{\comm}(\frac{P}{P_{\ge\infty}+I^\infty P}),\delta) @>\sim>>(\Omega (\frac{T}{F_\infty T+J^\infty}
),b)
$$
\endproclaim
\demo{Proof} Consider the associative product $x\star y =e^{-1}(exey)$ for $x,y\in SL=P$; put $Q:=(P,\star)$.
The map $e$ induces a chain isomorphism between $\Omega Q=(\Omega P,b^Q)$ and
$(\Omega T,b)$. Moreover, it follows from [Co3, 2.2], [C4 2.1] and [C4 2.6] that
the map $e$ induces a chain pro-isomorphism
$$
(\Omega(\frac{P}{P_{\infty}+I^\infty P}),b^Q)\cong (\Omega(\frac{T}{F_\infty T+J^\infty}),b)
$$
Let $\alpha:SdL=\Lambda L@>>> d\Omega P$, 
$$
\alpha(dg_1\land\dots\land dg_n)=\sum_{\sigma\in \Sigma_n}\sg(\sigma)dg_{\sigma 1}\dots dg_{\sigma n}
$$
Note $\eta=e\circ(1\otimes \alpha)$. By [Kas, Th. 3-a)], $1\otimes\alpha$ is a chain map. It induces a chain pro-map
$$
\Omega_{\comm}(\frac{P}{P_{\ge \infty}+I^\infty P})@>>>\Omega \frac{P}{P_{\ge\infty} +I^\infty P}\tag{61}
$$
We must show \thetag{61} is a homotopy equivalence. For this we shall construct a homotopy inverse 
$\beta:\Omega Q@>>>\Omega_{\comm}P$ of $1\otimes \alpha$ and homotopies $\gamma:\beta(1\otimes\alpha)@>>>1$ and
$\kappa:(1\otimes\alpha)\beta@>>>1$ each of which will be continous for the linear topologies of the filtrations
$\{\ker(\Omega_{\comm}P\fib\Omega_{\comm}P/P_{\ge n}+I^nP)\}_n$ and $\{\ker(\Omega P\fib\Omega P/P_{\ge n}+I^nP)\}_n$. We point
out that the first of these topologies coincides with that of the filtration 
$$\{I^n\Omega_{\comm}P+\bigoplus_{l\ge n}{}_l\Omega_{\comm}P\}_n$$
Write
$$
K: (\Omega_{\comm}\otimes P,\delta')@>\epsilon>> P
$$
for the augmented $Q\otimes Q^{op}$-resolution denoted $(L',b')$ in [Kas, Prop. 3] and
$$
R: (\Omega P\otimes P,b')@>\mu>>P
$$
for the augmented Hochschild resolution. By [Kas, Lemme 9] the continuous map $1\otimes \alpha\otimes 1:K@>>>R$
a chain $Q\otimes Q^{op}$-module homomorphism. It suffices to construct continuous $Q\otimes Q^{op}$-homomorphisms
$\beta':R@>>>K$, $\gamma':R@>>>R[1]$ and $\kappa':K@>>>K[1]$ such that $\beta'b'=\delta'\beta$, and such that
$\gamma'$ and $\kappa'$ be homotopies $\beta'(1\otimes\alpha\otimes 1)@>>>1$ and $1@>>>(1\otimes\alpha\otimes 1)\beta'$.
In turn for this it suffices to show that both $R$ and $K$ have continuous $k$-linear contracting homotopies. For then
the statndard procedure for lifting the identity in dimension zero to a chain map $\beta'$ using a contracting homotopy
for $K$ yields a continous $\beta'$, and similarly for the standard procedure for constructing the homotopies $\gamma'$
and $\kappa'$. The map $a\mapsto 1\otimes a$, $\omega\otimes x\mapsto dw\otimes x$ defines a continous contracting homotopy
for the augmented resolution $R$. To obtain a continous contracting homotopy for $K$ proceed as follows. Put
$$
_mK_*:=\bigoplus_{i+j=m}{}_i\Omega_{\comm}P\otimes P_j
$$
We have
$$
K=\bigoplus_{m=0}^\infty{}_mK\qquad \delta'(_mK)\subset \bigoplus_{p\ge m}{}_pK
$$
Let $\delta'_n$ be the homogeneous component of degree $n\ge 0$. By [Co3, 2.2] each $\delta'_n$ is a continous map.
Hence if $h_0$ is a continous contracting homotopy for 
$(K,\delta'_0)$ then the map $h=\sum_{m=0}^\infty h_m$ of the perturbation lemma [Co2, 2.5] is a continous contracting
homotopy for $(K,\partial')$. We remark that $(K,\partial'_0)$ is the standard Koszul resolution of $P$ as a module
over $P\otimes P$ with its commutative structure. Thus essentially the same argument as in the proof of lemma 7.5 gives
a continuous homotopy $h_0$ as wanted.\qed
\enddemo
\bigskip
\head{8. Proofs of the main theorems}\endhead
\bigskip
\demo{8.1.Proof of Theorem 6.2} 
We first do the case $l<\infty$. 
To start, we prove the existence of the isomorphism $\alpha_l$. Assume 
first $X=\Spec A$. Choose a presentation
$$
0@>>>J@>>>R@>\pi>>A@>>>0\tag{62}
$$
of $A$ as quotient of a quasi-free $R\in\Ass$. Put $R_l=R/F_{l+1}R$ and assume $\Cal Y$ is the 
following system of embeddings
$$
\Cal Y:=\{X\hookrightarrow Y:=\Spec R_l\}\tag{63}
$$
By lemma 7.2, for $f\in A$ and 
$R\supset \Gamma=\pi^{-1}\{f^n: n\ge 0\}$
we have canonical isomorphisms
$$
\Cal C_{\Cal Y}(\Cal O)^p(D(f))=\{\Cal O_{Y^p_n}(D(f))\}=C_{\ncl}(R[\Gamma^{-1}],J_f,\Cal O)
\tag{64}
$$
where $ J_f=\ker (R[\Gamma^{-1}]\fib A_f)$. Note $R[\Gamma^{-1}]$ is quasi-free, whence by proposition
7.1 we have a natural homotopy equivalence between the pro-complex \thetag{64} and
$$
\Omega_\ncl(R[\Gamma^{-1}],J_f)\tag{65}
$$
Now because the $D(f)$ form a basis for the Zariski topology, \thetag{65} determines a unique 
pro-complex of sheaves; by lemma 7.2 this pro-complex must be $\Omega^{\Cal Y}$. 
Next if 
$$
0@>>>J'@>>>R'_l@>>>A@>>>0
$$
is any presentation of $A$ as quotient of a formally $l$-smooth $R'_l\in\ncl$, then by [Co1, 3.3]
a choice of a map $R_l@>>>R'_l$ covering the identity of $A$ induces a homotopy equivalence
$$
\Omega_\ncl(R,J)@>\sim>>\Omega_\ncl(R',J')
$$
which in turn gives a homotopy equivalence  
$$
\Omega_{\ncl}^{\Cal Y}@>\sim>>\Omega_{\ncl}^{\Cal Y'}
$$
for $\Cal Y'=\{X\hookrightarrow Y'\}$. Now no longer assume $X$ is affine. If 
$\Cal Y$ consists entirely of affine embeddings, say  
$\Cal Y=\{\Spec A_i\hookrightarrow \Spec R^i_l\}$ for some affine open covering $\{\Spec A_i\}$ of
$\Spec A$, then each $n$-fold intersection 
$$
\Spec A_{i_0}\cap\dots\cap\Spec A_{i_n}\hookrightarrow
\Spec R_l^{i_j}\times_l\dots\times_l\Spec R_l^{i_n}=
\Spec (R_l^{i_0}*\dots *R_l^{i_n})_l
\tag{66}
$$
is of the form \thetag{63}, (by 2.4.3) so from the affine case we obtain 
a homotopy  equivalence
$$
\Cal C_{\Cal Y}(\Cal O)@>\alpha>\sim>\Omega^{\Cal Y}_{\ncl}
$$
If $\Cal Y$ is arbitrary then there is a finer system $\Cal Y'$ which consists entirely of affine
embeddings; the argument of [H1, Remark on page 28] shows the refinement map
$$
\Omega^{\Cal Y}_{\ncl}@>\alpha>\sim>\Omega^{\Cal Y'}_{\ncl}
$$ 
is a quasi-isomorphism. 
Next we construct the map $e$ of the theorem. Assume first $X=\Spec A$, choose a presentation
\thetag{62} with $R=TV$, a tensor algebra and let $\Cal Y$ be as in \thetag{63}.
Further consider the system of $\npl$-embeddings $\Cal Z$ consisting
of the single embedding $\Spec A\hookrightarrow \Spec \Poiss_{\le l}V$ induced by 
the composite
$$
\Poiss_{\le l}V=P_{\le l}SV\fib SV\cong \frac{R}{F_1R}\fib A
$$
Then by lemmas 7.2 and 7.3 and by [C4 2.1 and 2.6] the map \thetag{49} induces an isomorphism of 
pro-complexes of abelian sheaves
$$
\Cal C_{\Cal Z}(\calo)@>e>\sim>\Cal C_{\Cal Y}(\calo)
$$
This proves the isomorphism $e$ in the affine case. For general $X$ one chooses $\Cal Y$
to consist entirely of affine embeddings as above, and $\Cal Z$ as the associated $\ncl$- system;
then use the affine case and the argument of \thetag{66}, taking into account that the free
product of tensor algebras is again a tensor algebra. 
To construct the map $e'$ in the affine case one uses a DG-version of the same
argument as for the construction of $e$. As per the arguments above, the affine case generalizes to 
the case when $\Cal Y$ consists of affine $\ncl$-embeddings, and $\Cal Z$ is the associated 
$\npl$-system.
It has already been shown that the hypercohomology of $\Omega^{\Cal Y}_{\ncl}$ is independent of the 
choice of $\Cal Y$; the same argument shows that of $\Omega^{\Cal Z}_{\npl}$ is independent of
the choice of $\Cal Z$. To finish the proof it suffices to show that there is a choice of $\Cal Y$,
$\Cal Z$ and $\Cal W$ for which $\pi'\circ e'$ is a cohomology isomorphism. Choose an affine open covering
$\Cal U$ of $X$. For each $\Cal U\owns U=\Spec A_U$ choose a presentation $A_U=SV_U/I_U$ and
let 
$$
\align
\Cal W=&\{\Spec A_U\hookrightarrow\Spec SV_U:U\in\Cal U\}\\
\Cal Z=&\{\Spec A_U\hookrightarrow\Spec\Poiss_{\le l}V_U:U\in\Cal U\}\\
\Cal Y=&\{\Spec A_U\hookrightarrow\Spec T_{\le l}V_U:U\in\Cal U\}
\endalign
$$
Then by the argument of \thetag{66} we are reduced to showing that if $A=SV/I$ then the projection
$$
\multline
\theta:P_{\le l}\Omega_{\comm}(SV/I^\infty)\cong \frac{\Poiss_{\le l}(V\oplus dV)}
{I^{\infty}\Poiss_{\le l}(V\oplus dV)}
\fib \\
\frac{S(V\oplus dV)}{I^{\infty}S(V\oplus dV)}\cong \Omega_{\comm}(SV/I^{\infty})
\endmultline
$$
is a homotopy equivalence. This follows from lemma 7.5, and the fact that 
$$
\ker\theta=\bigoplus_{m=1}^l\frac{\Poiss_m(V\oplus dV)}{I^{\infty}\Poiss_m(V\oplus dV)}
$$
This concludes the proof of the case $l<\infty$ of the theorem. Because the cohomology isomorphisms we found come
from natural cochain equivalences which are compatible with the inclusions 
$\ncl\subset {NC}_{l+1}$
and $\npl\subset {NP}_{l+1}$, the case $l=\infty$ follows.\qed
\enddemo
\bigskip
\remark{8.2. Remark} A similar argument as that of the last part of the proof above shows that
$$
H^n(X_{\comm-\inf},P_m)=0\qquad (m\ge 1, n\ge 0)
$$
Indeed because $P_m$ is quasi-coherent it suffices --by \thetag{40}-- to show that 
for $A$ a commutative algebra,
$$
\align
\qquad\qquad H^n(\inf(\comm/A),P_m)=&0\tag{67}\\
\intertext{With the notations of the proof, we see using Prop. 7.1 that 
\thetag{67} equals}
=&H^n(\lim_r P_m\Omega_{\comm}(SV/I^r))\\
\endalign
$$
which is zero by lemma 7.5.
\endremark
\bigskip
\demo{8.3. Proof of Theorem 6.5} We shall assume $X$ affine and $l<\infty$.  The same argument
as in the proof of theorem 6.2 applies to deduce the theorem from this particular case. Let $X=\Spec A$; choose
$R$, $\pi$ and $J$ as in \thetag{62} and let $\Cal Y$ be as in \thetag{63}. Proceed as in the proof of
theorem 6.2 to obtain a natural homotopy equivalence of pro-complexes of vectorspaces
$$
\frac{\Cal C_{\Cal Y}(\calo)(D(f))}{[\Cal C_{\Cal Y}(\Cal O)(D(f)),\Cal C_{\Cal Y}(\Cal O)(D(f))]}
@>{\bar{\alpha}}_l
>\sim>\frac{\Omega_{\ncl}(R[\Gamma^{-1}],J_f)}{[\Omega_{\ncl}(R[\Gamma^{-1}],J_f),\Omega_{\ncl}(R[\Gamma^{-1}],J_f)]}
\tag{68}
$$
for each $f\in A$.
Because ${\bar{\alpha}}_l$ is natural and because sheafification depends only on the value of the presheaf on a basis of the 
topology,
\thetag{68} induces a homotopy equivalence of pro-complexes of sheaves
$$
\Cal C_{\Cal Y}(\frac{\calo}{[\calo,\calo]})=\frac{\Cal C_{\Cal Y}(\calo)}{[\Cal C_{\Cal Y}(\Cal O),\Cal C_{\Cal Y}(\Cal O)]}
@>{\bar{\alpha}}_l
>\sim>\frac{\Omega_{\ncl}^{\Cal Y}}{[\Omega_{\ncl}^{\Cal Y},\Omega_{\ncl}^{\Cal Y}]}
\tag{69}
$$
This cochain map gives the cohomology isomorphism of the theorem. The same argument as in the proof of theorem 6.2  shows
that the homotopy type of the complexes \thetag{69} is the same as that of those obtained from a different choice
of $\ncl$-embedding $\Cal Y'=\{X\hookrightarrow \Spec R'_l\}$. Assume now $R=TV$, a tensor algebra, and choose
$\Cal Y$, $\Cal Z$ and $\Cal W$ as in the proof of theorem 6.2. Then by lemmas 7.2, 7.3, and 7.4, 
[C4, 2.1 and 2.6] and sheafification, we have isomorphisms of pro-complexes of Zariski sheaves
$$
\bigoplus_{m=0}^l\Cal C_{\Cal W}((P_m)_\delta)
@<\gamma<\cong<\Cal C_{\Cal Z}(\frac{\calo}{\{\calo,\calo\}})
@>\bar{e}>\cong>\Cal C_{\Cal Y}(\frac{\calo}{[\calo,\calo]})
$$ 
where $P_m$ is as in remark 8.2 above and the subscript indicates the sheaf cokernel of the restriction 
of the 
coboundary map $\delta$ of 6.7 to $_{m-1}\Omega_{\comm}P$, i.e. the sheafification of
$$
D(f)\mapsto\frac{P_mA_f}{\sum_{i+j=m-1}\{P_iA_f,P_jA_f\}}.
$$
A DG-version of the same argument gives isomorphisms
$$
\bigoplus_{m=0}^l(P_m\Omega_{\comm}^{\Cal W})_\delta @<\gamma'<\cong<
\frac{\Omega_{\npl}^{\Cal Z}}{\{\Omega_{\npl}^{\Cal Z},\Omega_{\npl}^{\Cal Z}\}}@>{\bar{e}}'>\cong>
\frac{\Omega_{\ncl}^{\Cal Y}}{[\Omega_{\ncl}^{\Cal Y},\Omega_{\ncl}^{\Cal Y}]}
$$
By naturality we get a homotopy equivalence
$$
\Cal C_{\Cal W}((P_m)_\delta)@>>\sim>(P_m\Omega_{\comm}^{\Cal W})_\delta
$$
Next, consider the truncation $\tau_{2m}\N$ of the complex of Lemma 7.7. By 7.9, 8.2
and Proposition 7.9, we have a commutative diagram of homotopy equivalences 
$$
\CD
\op{Tot}{\Cal C}_{\Cal W}(\tau_m\N)@>>\sim>{\Cal C}_{\Cal W}((P_m)_\delta)\\
@VVV  @VVV\\
\op{Tot}{\Cal C}_{\Cal W}(\tau_m\Omega)@>>\sim>{\Cal C}_{\Cal W}(\frac{\Omega_{\comm}^m}{d\Omega_{\comm}^{m-1}})\\
\endCD
$$
with as rows the natural projections and as columns the maps induced by that of Lemma 7.7. To finish the proof we
must show that the natural projection
$$
\op{Tot}{\Cal C}_{\Cal W}(\tau_m\Omega_{\comm})
@>>\sim>{\Cal C}^0_{\Cal W}(\tau_m\Omega_{\comm})=\tau_m({\Omega_{\comm}}^{\Cal W})\tag{70}
$$
is a homotopy equivalence. Recall $\Cal W=\{\Spec A\hookrightarrow \Spec S\}$, where $S=SU$ is the symmetric algebra
of some vectorspace $U$ and $A=S/I$. Thus by \thetag{60} we have an isomorphism (where $Cyl$ is as in
[Co2])
$$
Cyl^p(SU)=\frac{S(U\oplus U\otimes V^p)}{< U\otimes V^p>^\infty}
$$
Grade $S(U\oplus U\otimes V^p)$ by $|(u,0)|=0$, $|(0,u\otimes v)|=1$. Then there is an inclusion
$$
\tau_m\Omega_{\comm}SU
\overset{\iota}\to{\hookrightarrow} \tau_m\Omega_{\comm} Cyl^*SU=C^*(SU,0,\tau_m\Omega_{\comm})
$$
of the constant co-simplicial cochain pro-complex as the part of degree zero of the \v Cech-Alexander pro-complex. 
The map $\iota$ is a right inverse for the natural projection 
$\mu:C^*(SU,0,\tau_m\Omega_{\comm})\fib\tau_m\Omega_{\comm}SU$. The Cartan homotopy associated to the degree derivation
$D(x)=|x|x$ gives a homotopy $\iota\mu@>>>1$ which 
is compatible with the cosimplicial structure, localization and the $I$-adic topology. Thus upon sheafification we
get that \thetag{70} is a homotopy equivalence.\qed
\enddemo
\bigskip
\remark{8.3.4}It follows from 7.1 and the proof above that if $R$ is formally $\ncl$-smooth ($l<\infty$)
and $A=R/F_{1}R$ then 
$$
\align
H^n(\inf(\ncl/A),\frac{\calo}{[\calo,\calo]})=
&H^n(\frac{\Omega_{\ncl}R}{[\Omega_{\ncl}R,\Omega_{\ncl}R]})\\
=&\bigoplus_{m=0}^nH_{dR}^{n+2m}A.
\endalign
$$
\endremark
\bigskip
\demo{8.4. Proof of Theorem 6.8}
Let $A\in\comm$, $X=\Spec A$, $R$, $J$ and $\pi$ as in 
\thetag{62}. For $l\ge 0$ put  $R_l=R/F_{l+1}R\supset J_l=J+F_{l+1}R/F_{l+1}R$, $Y_l=\Spec R_l$,
$\Cal Y=\{X\hookrightarrow Y_1\hookrightarrow Y_2\hookrightarrow\dots\}$ the associated formally $\nci$-smooth embedding.
By Goodwillie's theorem ([CQ2, Th. 10.1]), we have a natural isomorphism
$$
HC^{per}_*(X):=\Hy^*(X_{\pro-\zar},\Cal{CC}^{per})\cong\Hy^*(X_{\pro-\zar}, (\Cal {CC}^{per})^{\Cal Y}) 
$$
To prove the isomorphism $f_4$ it suffices to show that the natural pro-complex projection
$$
CC^{per}(R_\infty/J_\infty^\infty)\fib \fx(R_\infty/J_\infty^\infty)\tag{71}
$$
is a quasi-isomorphism. Because \thetag{71} comes from a map of mixed complexes, it will suffice to show that the
map between the corresponding Hochschild complexes
$$
(\Omega(R_\infty/J_\infty^\infty),b)\fib (\fx(R_\infty/J_\infty^\infty),b)
$$
is a quism. By [Co1, 3.3] it suffices to prove this in the case when $R=TV$, a tensor algebra. 
With the notation of Lemma 7.11, we have a commutative diagram
$$
\CD
(\Omega_{\comm}(PSV/I^\infty PSV+P_\infty SV),
\delta)@>\eta>\sim> (\Omega (R_\infty/J_\infty^\infty),b)\\
@VVV @VVV\\
\fy(PSV/I^\infty PSV+P_\infty SV),\delta)@>\eta>\cong>(\fx(R_\infty/J_\infty^\infty),b)\\
\endCD
$$
where the top and bottom rows are respectively a homotopy equivalence and an isomorphism by Lemma 7.11 and its proof.
By lemma 7.6, the first vertical arrow is a quism, whence so is the second. This gives isomorphisms
$f_4$ and $f_{3'}$. It follows that all the coface maps of the \v Cech-Alexander co-simplicial 
pro-complex $C_{\Cal Y}(\fx)$
are homotopy equivalences, which gives isomorphism $f_2$. A similar argument produces an
isomorphism $f_{2'}$ for 
$\Cal Z=\{X\hookrightarrow \Spec P_{\le 1}SV\hookrightarrow\Spec P_{\le 2}\hookrightarrow\dots\}$; the passage
from this to the case when $P_{\le\infty}SV$ is replaced by an arbitrary formally $\npi$-smooth pro-algebra is done as in 
the $\nci$-case. The Poisson grading \thetag{44} induces a pro-cochain complex decomposition
$$
\fy(\left\{\frac{PSV}{P_{\ge n}SV}\right\}_n)\cong\frac{\fy(PSV)}{F_\infty\fy(PSV)}=\{\prod_{l=0}^n\N [-2l]\}_n
\tag{72}
$$
From \thetag{72}, lemma 7.7 and \thetag{33} we obtain the isomorphism $f_1$. 
This finishes the proof of the theorem
in the affine case; the general case follows from this by the same argument as in theorems 6.2 and 6.5.
\qed\enddemo

\enddocument